\documentclass[12pt,twoside]{article} 
\usepackage[english]{babel}
\usepackage[latin1]{inputenc} 
\usepackage{amsmath}
\usepackage{amssymb,amsfonts}
\usepackage{graphicx}                   
\usepackage{times,amssymb,amscd} 
\newcommand{\bA}{\mathbf{A}}

\newcommand{\bE}{\mathbf{E}}
\newcommand{\bG}{\mathbf{G}}
\newcommand{\bH}{\mathbf{H}}

\newcommand{\bL}{\mathbf{L}}
\newcommand{\bN}{\mathbf{N}}
\newcommand{\bR}{\mathbf{R}}			     
\newcommand{\bS}{\mathbf{S}}
\newcommand{\bV}{\mathbf{V}}

\newcommand{\be}{\mathbf{e}}

\newcommand{\bx}{\mathbf{x}}
\newcommand{\by}{\mathbf{y}}

\newcommand{\EUC}{\bE^3}
\newcommand{\SPH}{\bS^3}
\newcommand{\HYP}{\bH^3}
\newcommand{\SXR}{\bS^2\!\times\!\bR}
\newcommand{\HXR}{\bH^2\!\times\!\bR}
\newcommand{\SLR}{\widetilde{\bS\bL_2\bR}}
\newcommand{\NIL}{\mathbf{Nil}}
\newcommand{\SOL}{\mathbf{Sol}}

\begin{document}

\pagestyle{myheadings}
\markboth{\centerline{Emil Moln\'ar and Jen\H o Szirmai}}
{Classification of $\SOL$ lattices}
\title
{Classification of $\SOL$ lattices}

\author{\normalsize{Dedicated to our Friend,} \\
\normalsize{\medbreak Professor Erhard QUAISSER 
on the Occasion of His 70th Birthday} \\ \medbreak Emil Moln\'ar and Jen\H o Szirmai\thanks{The research of the author is supported by Action \"Osterreich Ungarn 2008 No.71öu1,  
and by the Slovenian -- Hungarian Intergovernmental Cooperation Programme 2008-2009.} \\
\normalsize Budapest University of Technology and \\ 
\normalsize Economics Institute of Mathematics, \\
\normalsize Department of Geometry \\
\normalsize H-1521 Budapest, Hungary \\
\date{\normalsize \today}}

\maketitle

\begin{abstract}
$\SOL$ geometry is one of the eight homogeneous Thurston 3-geomet\-ri\-es 
$$\EUC, \SPH, \HYP, \SXR, \HXR, \SLR, \NIL, \SOL.$$
In \cite{Sz10} the {\it densest lattice-like translation ball packings} to a type (type {\bf I/1} in this paper) of $\SOL$ lattices has been determined.
Some basic concept of $\SOL$ were defined by {\sc{P. Scott}} in \cite{S}, in general.  

In our present work we shall classify $\SOL$ lattices in an algorithmic way into 17 (seventeen) types, 
in analogy of the 14 Bravais types of the Euclidean 3-lattices, but infinitely many $\SOL$ affine 
equivalence classes, in each type.
Then the discrete isometry groups of compact fundamental domain (crystallographic groups) can also be classified into infinitely 
many classes but finitely many types, left to other publication.
To this we shall study relations between 
$\SOL$ lattices and lattices of the pseudoeuclidean (or here rather called Minkowskian) plane \cite{AQ}. 
Moreover, we introduce the notion of $\SOL$ parallelepiped to every lattice type. From our new results we emphasize Theorems 3-4-5-6. 

In this paper we shall use the affine model of $\SOL$ space through affine-projective homogeneous coordinates \cite{M97} which 
gives a unified way of investigating and visualizing homogeneous spaces, in general. 
\end{abstract}
\newtheorem{theorem}{Theorem}[section]
\newtheorem{corollary}[theorem]{Corollary}
\newtheorem{lemma}[theorem]{Lemma}
\newtheorem{exmple}[theorem]{Example}
\newtheorem{defn}[theorem]{Definition}
\newtheorem{rmrk}[theorem]{Remark}
\newtheorem{proposition}[theorem]{Proposition}
\newenvironment{definition}{\begin{defn}\normalfont}{\end{defn}}
\newenvironment{remark}{\begin{rmrk}\normalfont}{\end{rmrk}}
\newenvironment{example}{\begin{exmple}\normalfont}{\end{exmple}}
\newenvironment{acknowledgement}{Acknowledgement}
\section{On Sol geometry}
\label{sec:1}
In this Section 1 we summarize the significant notions and notations of real $\SOL$ geometry (see \cite{M97}, \cite{S}).

$\SOL$ is defined as a 3-dimensional real Lie group with multiplication  
\begin{equation}
     \begin{gathered}
(a,b,c)(x,y,z)=(x + a e^{-z},y + b e^z ,z + c).
     \end{gathered} \tag{1.1}
     \end{equation}
We note that the conjugacy by $(x,y,z)$ leaves invariant the plane $(a,b,c)$ with fixed $c$:
\begin{equation}
     \begin{gathered}
(x,y,z)^{-1}(a,b,c)(x,y,z)=(x(1-e^{-c})+a e^{-z},y(1-e^c)+b e^z ,c).
     \end{gathered} \tag{1.2}
     \end{equation}
Moreover, for $c=0$, the action of $(x,y,z)$ is only by its $z$-component,where $(x,y,z)^{-1}=(-x e^{z}, -y e^{-z} ,-z)$. Thus the $(a,b,0)$ plane is distinguished as a {\it base plane} in
$\SOL$, or by other words, $(x,y,0)$ is normal subgroup of $\SOL$.
$\SOL$ multiplication can also be affinely (projectively) interpreted by "right translations" 
on its points as the following matrix formula shows, according to (1.1):
     \begin{equation}
     \begin{gathered}
     (1;a,b,c) \to (1;a,b,c)
     \begin{pmatrix}
         1&x&y&z \\
         0&e^{-z}&0&0 \\
         0&0&e^z&0 \\
         0&0&0&1 \\
       \end{pmatrix}
       =(1;x + a e^{-z},y + b e^z ,z + c)
       \end{gathered} \tag{1.3}
     \end{equation}
by row-column multiplication.
This defines "translations" $\mathbf{L}(\mathbf{R})= \{(x,y,z): x,~y,~z\in \mathbf{R} \}$ 
on the points of space $\SOL= \{(a,b,c):a,~b,~c \in \mathbf{R}\}$. 
These translations are not commutative, in general. 
Here we can consider $\mathbf{L}$ as projective collineation group with right actions in homogeneous 
coordinates as usual in classical affine-projective geometry.
We will use the Cartesian homogeneous coordinate simplex $E_0(\be_0)$,$E_1^{\infty}(\be_1)$, \ $E_2^{\infty}(\be_2)$, \ 
$E_3^{\infty}(\be_3), \ (\{\be_i\}\subset \bV^4$ \ $\text{with the unit point}$ $E(\be = \be_0 + \be_1 + \be_2 + \be_3 ))$ 
which is distinguished by an origin $E_0$ and by the ideal points of coordinate axes, respectively. 
Thus {$\SOL$} can be visualized in the affine 3-space $\bA^3$
(so in Euclidean space $\bE^3$) as well \cite{BSz}.

In this affine-projective context E. Moln\'ar has derived in \cite{M97} the usual infinitezimal arc-length square at any point 
of $\SOL$, by pull back translation, as follows
\begin{equation}
   \begin{gathered}
      (ds)^2:=e^{2z}(dx)^2 +e^{-2z}(dy)^2 +(dz)^2.
       \end{gathered} \tag{1.4}
     \end{equation}

Hence we get infinitezimal Riemann metric invariant under translations, by the symmetric metric tensor field $g$ on $\SOL$ by components as usual.

It will be important for us that the full isometry group Isom$(\SOL)$ has eight components, since the stabilizer of the origin 
is isomorphic to the dihedral group $\mathbf{D_4}$, generated by two involutive (involutory) transformations, preserving (1.4): 
\begin{equation}
   \begin{gathered}
      (1)  \ \ y \leftrightarrow -y; \ \ (2)  \ x \leftrightarrow y; \ \ z \leftrightarrow -z; \ \ \text{i.e. first by $3\times 3$ matrices}:\\      
     (1) \ \begin{pmatrix}
               1&0&0 \\
               0&-1&0 \\
               0&0&1 \\
     \end{pmatrix}; \ \ \ 
     (2) \ \begin{pmatrix}
               0&1&0 \\
               1&0&0 \\
               0&0&-1 \\
     \end{pmatrix}; \\
     \end{gathered} \tag{1.5}
     \end{equation}
     with its product, generating a cyclic group $\mathbf{C_4}$ of order 4
     \begin{equation}
     \begin{gathered}
     \begin{pmatrix}
                    0&1&0 \\
                    -1&0&0 \\
                    0&0&-1 \\
     \end{pmatrix};\ \ 
     \begin{pmatrix}
               -1&0&0 \\
               0&-1&0 \\
               0&0&1 \\
     \end{pmatrix}; \ \ 
     \begin{pmatrix}
               0&-1&0 \\
               1&0&0 \\
               0&0&-1 \\
     \end{pmatrix};\ \ 
     \mathbf{Id}=\begin{pmatrix}
               1&0&0 \\
               0&1&0 \\
               0&0&1 \\
     \end{pmatrix}. 
     \end{gathered} \notag
     \end{equation}
     Or we write by collineations fixing the origin $O(1,0,0,0)$:
\begin{equation}
(1) \ \begin{pmatrix}
         1&0&0&0 \\
         0&1&0&0 \\
         0&0&-1&0 \\
         0&0&0&1 \\
       \end{pmatrix}, \ \
(2) \ \begin{pmatrix}
         1&0&0&0 \\
         0&0&1&0 \\
         0&1&0&0 \\
         0&0&0&-1 \\
       \end{pmatrix} \ \ \text{of form (1.3)}. \tag{1.6}       
\end{equation}
A general isometry of $\SOL$ to the origin $O$ is defined by a product $\gamma_O \tau_X$, first $\gamma_O$ of form (1.6) then $\tau_X$ of (1.3). To
a general point $A(1,a,b,c)$, this will be a product $\tau_A^{-1} \gamma_O \tau_X$, mapping $A$ into $X(1,x,y,z)$. 

Conjugacy of translation $\tau$ by an above isometry $\gamma$, as $\tau^{\gamma}=\gamma^{-1}\tau\gamma$ also denotes it, will also be used by 
(1.3) and (1.6) or also by coordinates with above conventions.

We remark only that the role of $x$ and $y$ can be exchanged throughout the paper, but this leads to mirror interpretation of $\SOL$.
As formula (1.4) fixes the metric of $\SOL$, the change above is not an isometry of a fixed $\SOL$ interpretation. Other conventions are also accepted
and used in the literature.
\section{Discrete subgroups of $Isom(\SOL)$}
\label{sec:2}
\begin{definition}
Let $\Gamma < \mathbf{L(R)}$ be a subgroup, generated by three independent translations $\tau_1(t_1^1,t_1^2,t_1^3)$,
$\tau_2(t_2^1,t_2^2,t_2^3)$, $\tau_3(t_3^1,t_3^2,t_3^3)$ with non-commutative addition, or in this sense with $\mathbf{Z}$ (integer) linear combinations. 
Here the upper indices indicate the corresponding $(\mathbf{e}_1,\mathbf{e}_2,\mathbf{e}_3)$ coordinates of basis translations.
$\Gamma$ is called discrete translation group or lattice of $\SOL$, if its action is discrete (by the induced orbit topology), i.e.
there is a compact fundamental "parallelepiped" (with side face pairing or identifications of its "bent" side faces) 
$\widetilde{\mathcal{F}}=\SOL/ \Gamma$
(Fig.~2), representing a compact $\SOL$ manifold $M$. $\square$ 
\end{definition}
More generally, we have
\begin{definition}
Let $\mathbf{G} < Isom(\SOL)$ be a subgroup of isometries of $\SOL$, such that $\mathbf{G}$ has a compact fundamental domain 
$M \sim \widetilde{\mathcal{F}}_{\mathbf{G}} \sim \SOL/\mathbf{G}$ with face pairing (face identifications) which represents the orbit space 
$M$ of $\SOL$ under $\mathbf{G}$. $\square$
\end{definition}
As in the classical Euclidean geometry, we have the analogue of Bieberbach theorem also in $\SOL$ geometry, i.e. 
above $\mathbf{G}$ has a 3- dimensinal translation subgroup, i.e. a 3-lattice in $\SOL$.
In \cite{S} {\sc P. Scott \cite{S} } indicated the proof of  
\begin{theorem}[P. Scott] 
Let $\bG$ be a discrete subgroup of $Isom(\SOL)$ which acts freely on $\SOL$ with quotient $M$. Then the natural 2-dimensional foliation of 
$M$ gives $M$ the structure of 2-dimensional bundle over a 1-dimensional orbifold. $\square$
\end{theorem}
We have the more descriptive two theorems 2-3 for translation lattices:

\begin{theorem}[\cite{S}, \cite{B}, \cite{Sz10}] 
Each lattice $\Gamma$ of $\SOL$ has a group presentation 
\begin{equation}
\begin{gathered}
\Gamma=\Gamma(\Phi)=\{\tau_1,\tau_2,\tau_3: [\tau_1,\tau_2]=1, \tau_3^{-1}\tau_1\tau_3=\tau_1 \Phi^T, \tau_3^{-1}\tau_2\tau_3=\tau_2 \Phi^T
\}, \tag{2.1}
\end{gathered}
\end{equation}
where 
\begin{equation}
\Phi= \begin{pmatrix}
         p & q \\
         r & s \\
\end{pmatrix}\in SL_2(\mathbf{Z}) \ with \ \mathrm{tr}(\Phi)=N=p+s > 2, \ ps-qr=1, \tag{2.2}
\end{equation}
such that for  
$\tau_1(t_1^1,t_1^2,t_1^3), \ \tau_2(t_2^1,t_2^2,t_2^3)$ above the equalities $t_1^3=0=t_2^3$ hold, and 
the matrix  
\begin{equation}
T=\begin{pmatrix}
         t_1^1 & t_1^2 \\
         t_2^1 & t_2^2 \\
\end{pmatrix} \in GL_2(\mathbf{R}) \ \ satisfies: \ \ 
T^{-1}\Phi T =:\Phi^T= 
\begin{pmatrix}
         e^{-t_3^3}&0 \\
         0& e^{t_3^3} \\
\end{pmatrix}, \tag{2.3}
\end{equation}
that is just a hyperbolic rotation fixed by the component $t_3^3$ in $\tau_3$ above. $\square$
\end{theorem}
Translations $\tau_1$ and $\tau_2$ satisfy by their $x$ and $y$ coordinates   
the following equations $\tau_1 \Phi^T \sim (t_1^1 e^{-t_3^3}, t_1^2 e^{t_3^3})
\begin{pmatrix}
         \be_1 \\
         \be_2 \\
\end{pmatrix},$
$\tau_2 \Phi^T \sim (t_2^1 e^{-t_3^3}, t_2^2 e^{t_3^3})
\begin{pmatrix}
         \be_1 \\
         \be_2 \\
\end{pmatrix}$ 
in the commutative base vector plane of $\SOL$, spanned by $\be_1$ and $\be_2$. Here $\sim$ refers, for a while as a hint, to 
the distinction of affine coordinates. Later on this distinction will be clear (hopefully) from the context.
The basis vectors $\be_1$ and $\be_2$ are just the eigenvectors of $\Phi^T$ to eigenvalues
$e^{-t_3^3}$ and $e^{t_3^3}$, respectively. 
Our main tool in this paper will be
\begin{theorem} 
Translations $\tau_1, \tau_2, \tau_3$ form a lattice in $\SOL$ space if and only if their above affine coordinates satisfy the conditions summarized 
in the following two main cases:
\begin{enumerate}
\item[I.] $t_1^3=t_2^3=t_3^1=t_3^2=0$. (The so-called {\it "fundamental lattices"}, with above $p=0,~s=N,~q=1,~r=-1$ (see \cite{AQ}, \cite{Sz10}), 
form a special class in this case.)
\item[II.] $t_1^3=t_2^3=0$ \ and $(t_3^1,~t_3^2) \in \mathbf{R}^2 \setminus \{(0,0)\}$, taken $\text{mod}\langle \tau_1,\tau_2 \rangle$, i.e.
modulo by the invariant base plane lattice $\Gamma^0(\Phi)$.
\end{enumerate}
In the above cases $t_3^3=\log{\frac{N+\sqrt{N^2-4}}{2}}$ depends on the main parameter $\mathrm{tr}(\Phi)=N=p+s > 2$, introduced in Theorem 2.
\end{theorem}
{\it Proof:}

We consider three independent $\SOL$ translations, first in general form $\tau_1(t_1^1,t_1^2,t_1^3)$,
$\tau_2(t_2^1,t_2^2,t_2^3)$, $\tau_3(t_3^1,t_3^2,t_3^3)$ defined in (1.1) and (1.3). 

{\it We take the commutator translation of} 
$\tau_1, \tau_2$ by (1.3). 
\begin{equation}
\begin{gathered}
\lbrack \tau_1,\tau_2 \rbrack = \tau_1^{-1} \tau_2^{-1} \tau_1 \tau_2 = \\
=\begin{pmatrix}
1& & &t_2^1(1-e^{-t_1^3})+t_1^1(e^{-t_2^3}-1)& & t_2^2(1-e^{t_1^3})+t_1^2(e^{t_2^3}-1)& & &0 \\
0& & &1 & & 0 & & &0 \\
0& & &0 & & 1 & & &0 \\
0& & &0 & & 0 & & &1 \\
\end{pmatrix}. \tag{2.4}
\end{gathered}
\end{equation} 
If e.g. $t_1^3 \ne 0$ and $(t_2^1,t_2^2) \ne (0,0)$, this is a non-zero translation in the base plane $\langle \be_1,\be_2 \rangle$. 
Thus, $t_1^3=0$ can be and will be assumed in the following. Similarly,
the commutator $[\tau_1,\tau_3]$ lies also in the base plane, and is not trivial if $t_3^3 \ne 0$ and $(t_1^1,t_1^2) \ne (0,0)$. Moreover
in the last case, if $t_1^1 \ne 0 \ne t_1^2$ hold in addition, we have 
$[\tau_1,\tau_3] \sim t_1^1(e^{-t_3^3}-1) \be_1+t_1^2(e^{t_3^3}-1) \be_2$ as independent translation from $\tau_1$.
Hence we can consider a basis $\tau_1(t_1^1,t_1^2,0), \tau_2(t_2^1,t_2^2,0), \tau_3(t_3^1,t_3^2,t_3^3)$ for our lattice. 
{\it Then $[\tau_1,\tau_2]$ is zero translation, and
$t_3^1\be_1+t_3^2 \be_2$, can be assumed $\mod$ base lattice, generated by} $(\tau_1,\tau_2)$. 

{\it We take the second condition} of discrete action of $\Gamma$, requirements 
$\tau_3^{-1}\tau_1\tau_3=\tau_1 \Phi^T,$ $\tau_3^{-1}\tau_2\tau_3=\tau_2 \Phi^T$
(see 2.1). 
That means, translations $\tau_1(t_1^1,t_1^2,0)$, $\tau_2(t_2^1,t_2^2,0)$ generate a plane lattice
$\Gamma^0$ mapping the plane $[x,y]$  onto itself (given by its orthonormal coordinate system  $(O,\mathbf{e}_1,\mathbf{e}_2)$, 
see formula (1.4) and
Fig.~1). The third basis translation $\tau_3$ has the crutial component $t_3^3$ to the third $\mathbf{e}_3$ direction. 
The coordinate plane $(O,x,y)$ and the
parallel planes $( (O,x,y)^{(\tau_3)^k}, \ k \in \mathbf{Z})$ produced by $\tau_3$, 
contain congruent integer lattices generated by $\tau_1$ and $\tau_2$.
These point lattices are denoted also by $\Gamma^{k t_3^3} \ (k \in \mathbf{Z})$. 
We repeat for explanations the following formulas with computations: 
\begin{equation}
\begin{gathered} \begin{pmatrix} \tau_1 \\ \tau_2 \\ \end{pmatrix} \Phi^T= \begin{pmatrix} t_1^1&t_1^2 \\ t_2^1&t_2^2 \\ \end{pmatrix} 
\begin{pmatrix}
\be_1 \\ \be_2 \\ \end{pmatrix} 
\Phi^T= \begin{pmatrix} t_1^1&t_1^2 \\ t_2^1&t_2^2 \\ \end{pmatrix} \begin{pmatrix} e^{-t_3^3}&0 \\ 0&e^{t_3^3} \\
\end{pmatrix} \begin{pmatrix} \be_1 \\ \be_2 \\ \end{pmatrix} \in \Gamma^0(\tau_1,\tau_2), \end{gathered} \notag \end{equation} 
i.e. 
\begin{equation}
\begin{gathered} \begin{pmatrix} \tau_1 \\ \tau_2 \\ \end{pmatrix} \Phi^T= \begin{pmatrix} p&q \\ r&s \\ \end{pmatrix} \begin{pmatrix} \tau_1 \\
\tau_2 \\ \end{pmatrix} = \begin{pmatrix} p&q \\ r&s \\ \end{pmatrix} \begin{pmatrix} t_1^1&t_1^2 \\ t_2^1&t_2^2 \\ \end{pmatrix} \begin{pmatrix} \be_1 \\
\be_2 \\ \end{pmatrix}. \end{gathered} \tag{2.5} \end{equation}
By (1.3) we get
\begin{equation} 
\begin{gathered} \tau_3^{-1}\tau_1\tau_3 =
\begin{pmatrix} 1 & & t_1^1 e^{-t_3^3} & & t_1^2 e^{t_3^3} & & 0 \\ 0 & & 1& & 0 & & 0 \\ 0 & & 0 & & 1 & & 0 \\ 0 & & 0 & & 0 & & 1 \\ \end{pmatrix}=
\tau_1 \Phi^T, \\
\tau_3^{-1}\tau_2\tau_3
=\begin{pmatrix}
1 & & t_2^1 e^{-t_3^3} & & t_2^2 e^{t_3^3}& & 0 \\
0 & & 1 & & 0 & & 0 \\
0 & & 0 & & 1 & & 0 \\
0 & & 0 & & 0 & & 1 \\
\end{pmatrix}= \tau_2 \Phi^T. \tag{2.6}
\end{gathered} 
\end{equation}

Thus we obtain our "discreteness conditions" (2.2) and (2.3) for the coordinates of translations $\tau_1,~ \tau_2,~ \tau_3$ by equations 
(2.4), (2.5) and (2.6). 

Relations between parameters $p,q,r,s$ and $t_1^1,t_1^2,t_2^1,t_2^2$ will be shown in the next equations:
\begin{equation}
\begin{gathered}
\Phi=\begin{pmatrix}
p&q \\
r&s \\
\end{pmatrix}\in SL_2(\mathbf{Z}), \ \ \mathrm{i.e.} \ \ p,q,r,s \in \mathbf{Z}, \ \ \mathrm{and} \ \ ps-qr=1,  \ \text{so that} \\
0 \ne \det \begin{pmatrix}
t_1^1&t_1^2 \\
t_2^1&t_2^2 \\
\end{pmatrix}=t_1^1 t_2^2-t_1^2 t_2^1=D \ \Rightarrow 
\begin{pmatrix}
t_1^1&t_1^2 \\
t_2^1&t_2^2 \\
\end{pmatrix}^{-1}=
\frac{1}{D} \cdot \begin{pmatrix}
t_2^2&-t_1^2 \\
-t_2^1&t_1^1 \\
\end{pmatrix}. \ 
\end{gathered} \notag 
\end{equation}
Therefore, we get 
\begin{equation}
\begin{gathered}
\begin{pmatrix}
p&q \\
r&s \\
\end{pmatrix}=\frac{1}{D} \cdot 
\begin{pmatrix}
t_1^1&t_1^2 \\
t_2^1&t_2^2 \\
\end{pmatrix}
\begin{pmatrix}
e^{-t_3^3}&0 \\
0&e^{t_3^3} \\
\end{pmatrix} \tag{2.7}
\begin{pmatrix}
t_2^2&-t_1^2 \\
-t_2^1&t_1^1 \\
\end{pmatrix}=\\
=\frac{1}{D} \cdot \begin{pmatrix}
t_1^1 e^{-t_3^3} t_2^2-t_1^2 e^{t_3^3}t_2^1& &-t_1^1 e^{-t_3^3} t_1^2+t_1^2 e^{t_3^3}t_1^1 \\
t_2^1 e^{-t_3^3} t_2^2-t_2^2 e^{t_3^3}t_2^1& &-t_2^1 e^{-t_3^3} t_1^2+t_2^2 e^{t_3^3}t_1^1 \\
\end{pmatrix}= \\
=\begin{pmatrix}
e^{-t_3^3}-\frac{2t_1^2 t_2^1}{D} \sinh{t_3^3} & &\frac{2 t_1^1 t_1^2}{D} \sinh{t_3^3} \\
\frac{-2 t_2^1 t_2^2}{D} \sinh{t_3^3} & &e^{t_3^3}+\frac{2 t_1^2 t_2^1}{D} \sinh{t_3^3} \\
\end{pmatrix}.\\
\end{gathered} 
\end{equation}
From (2.7) we get $N=p+s=e^{-t_3^3}+e^{t_3^3}=2\cosh{t_3^3}$, then 
\begin{equation}
\frac{N-2p-\sqrt{N^2-4}}{2\sqrt{N^2-4}}=\frac{1}{\frac{t_1^1 t_2^2}{t_1^2 t_2^1}-1}, \ \ 
\frac{q}{\sqrt{N^2-4}}=\frac{1}{\frac{t_2^2}{t_1^2}-\frac{t_2^1}{t_1^1}} \tag{2.8}
\end{equation}
can be read off. Parameters $t_1^1,t_1^2,t_2^1,t_2^2$ can be expressed up to some freedom. Namely, we obtain
\begin{equation}
\frac{t_2^1}{t_1^1}=\frac{N-2p-\sqrt{N^2-4}}{2q}, \ \ \frac{t_2^2}{t_1^2}=\frac{N-2p+\sqrt{N^2-4}}{2q}, \tag{2.9}
\end{equation}
where $t_1^1 t_1^2=\frac{Dq}{\sqrt{N^2-4}}$ with $D=t_1^1t_2^2-t_1^2t_2^1$ show the freedom. Take e.g. 
$\tau_1=t_1^1\be_1+t_1^2\be_2=\frac{\sqrt{2}}{2}\be_1+\frac{\sqrt{2}}{2}\be_2=\be_{1}^,$. 

Then $\tau_2=t_2^1\be_1+t_2^2\be_2=\frac{1}{2q}(N-2p-\sqrt{N^2-4})\frac{\sqrt{2}}{2}\be_1+\frac{1}{2q}(N-2p+\sqrt{N^2-4})\frac{\sqrt{2}}{2}\be_2$
provides by $\be_{2}^,=-\frac{\sqrt{2}}{2}\be_1+\frac{\sqrt{2}}{2}\be_2$ (Fig.~1) 
\begin{equation}
\tau_2= \frac{1}{2q}(N-2p)\be_{1}^,+\frac{1}{2q}\sqrt{N^2-4}\be_{2}^, \tag{2.9'}
\end{equation}
with lattice parameters $N,~p,~q$ to be discussed later at (2.11). $\square$
\subsection{On Minkowskian plane lattices} 
If we ask for a $\Phi$-invariant bilinear form or scalar product
\begin{equation}
\langle \bx,\by \rangle = \langle x^i \be_i, y^j \be_j \rangle = x^i \langle \be_i, \be_j \rangle y^j=:x^i b_{ij} y^j \notag
\end{equation} 
(by Einstein-Schouten index conventions for indices 1, 2) for the typical lattice transformation $\Phi^T$ as a hyperbolic rotation, then we get the signature $(-, +)$ for $b_{ij}$. 
Thus, we say that the plane lattice $\Gamma^0(\Phi)$ is a {\it pseudoeuclidean or Minkowskian lattice}. In a new basis 
\begin{equation}
\begin{gathered}
\begin{pmatrix}
\be_1^,\\
\be_2^,\\
\end{pmatrix}=
\begin{pmatrix}
\frac{\sqrt{2}}{2}& \frac{\sqrt{2}}{2} \\
-\frac{\sqrt{2}}{2}& \frac{\sqrt{2}}{2} \\
\end{pmatrix}
\begin{pmatrix}
\be_1 \\
\be_2 \\
\end{pmatrix},
\end{gathered} \notag
\end{equation}
we get
\begin{equation}
\begin{gathered}
\begin{pmatrix}
\be_1^,\\
\be_2^,\\
\end{pmatrix} \Phi^{T'} =
\begin{pmatrix}
\cosh{t_3^3} & \sinh{t_3^3}\\
\sinh{t_3^3} & \cosh{t_3^3}\\
\end{pmatrix}
\begin{pmatrix}
\be_1^, \\
\be_2^, \\
\end{pmatrix} \ \ \text{for the hyperbolic rotation} \\
\mathrm{instead \ of} \ \ 
\begin{pmatrix}
\be_1 \\
\be_2 \\
\end{pmatrix} \Phi^{T}=
\begin{pmatrix}
e^{-t_3^3}& 0 \\
0 & e^{t_3^3} \\
\end{pmatrix}
\begin{pmatrix}
\be_1 \\
\be_2 \\
\end{pmatrix}.
\end{gathered} \tag{2.10}
\end{equation}
\begin{figure*}[]
\includegraphics[width=0.9\textwidth]{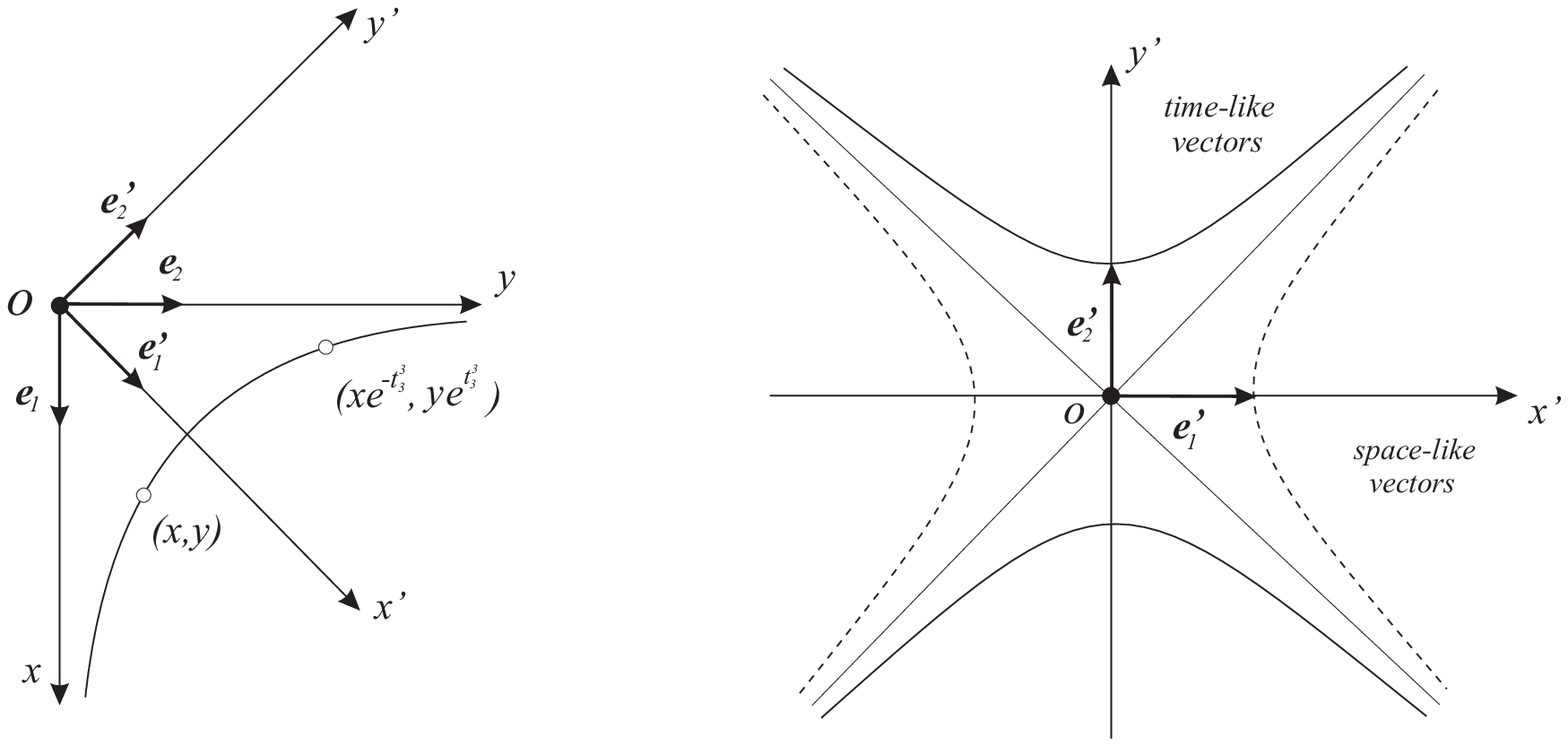}
\caption{}
\label{Fig:1}
\end{figure*}
The so-called {\it regular lattices} of pseudoeuclidean or Minkowskian plane were discussed by {\sc K. Alpers} and {\sc E. Quaisser} in \cite{AQ}
(see also \cite{BG}). 
Pseudoeuclidean plane, denoted by $\mathbf{E}_1^2$, is an affine plane over the field $\mathbf{R}$ of real numbers together with a symmetric bilinear form just
of signature $(-,+)$.
In \cite{AQ} these $\mathbf{E}_1^2$ plane lattices were classified, by the following:
\begin{definition}
A lattice $\Gamma_1$ in the Minkowskian plane is called affinely equivalent to a lattice $\Gamma_2$ $(\Gamma_1 \sim \Gamma_2)$, if there is an affine transformation $\alpha$ with
$\Gamma_2=\Gamma_1 \alpha$, so that any symmetry $S^\alpha$ of $\Gamma_2$ is derived by a symmetry $S$ of $\Gamma_1$ with conjugacy 
$S^ \alpha = \alpha^{-1} S \alpha$. A plane lattice $\Gamma^0$ is called {\it regular} if there is a hyperbolic rotation $(\Phi \sim \Phi^T=T^{-1} \Phi T$
as symmetry of $\Gamma^0$. $\square$
\end{definition}

In the sense of {\sc K. Alpers} and {\sc E. Quaisser} \cite{AQ}, with slight modification here, we are interested now in $\SOL$ lattices under the full isometry group
of $\SOL$ by (1.5) and (1.6) where its commutative plane lattices 
$\Gamma^0$ (generated by $\tau_1$ and $\tau_2$) are 
{\it regular}, i.e. admit hyperbolic rotations different from the identity(see (2.10)). 
The regular affine {\it lattice classes} $\Gamma^0(\Phi)(\tau_1,\tau_2)$ in the Minkowskian plane can be given by Definition 3 and (2.9')  
in the orthonormal coordinate system $\{O, \mathbf{e}_1^,, \mathbf{e}_2^, \}$ :
\begin{equation}
\begin{gathered}
\tau_1 \sim (1,0), \ \ \tau_2 \sim \frac{1}{2q}(N-2p,\sqrt{N^2-4}). \tag{2.11}
\end{gathered}
\end{equation}
\begin{definition}
An above lattice $\Gamma^0(\Phi)(\tau_1,\tau_2)$ which fulfils the
conditions: 
\begin{enumerate}
\item $p,q,r,s \in \mathbf{Z}$, are integers where $ps-qr=1$, moreover $p+s=:N=2\cosh{t_3^3} \geq 3$;
\item $0\le p \le \lfloor \frac{N}{2} \rfloor$ and $0<q$;    
\end{enumerate}
is called a $(N,p,q)$ lattice for a given trace $N=p+s$. 
\end{definition}
Of course, such a $(N,p,q)$ lattice with fixed $N=p+s$ and its affine class by Definition 3 are classically determined up to an unimodular, 
i.e. $SL_2(\mathbf{Z})$ conjugacy since a lattice basis is mapped to a lattice basis by any affinity, then $SL_2(\mathbf{Z})$ conjugacy
expresses any symmetry of $\Gamma^0$. Thus by a basis change: 
\begin{equation}
\begin{gathered}
\begin{pmatrix}
p'&q' \\
r'&s' \\
\end{pmatrix}=
\begin{pmatrix}
u&v \\
w& \overline{w} \\
\end{pmatrix}^{-1}
\begin{pmatrix}
p&q \\
r&s \\
\end{pmatrix}
\begin{pmatrix}
u&v \\
w& \overline{w} \\
\end{pmatrix}= \\
=\begin{pmatrix}
p \pm[v(N-2p)w-urv+wq\overline{w}] & \pm[ v(2p-N)\overline{w}-vrv+\overline{w}q\overline{w}] \\
\pm[u(N-2p)w+uru-wqw] & N-p \pm [v(2p-N)w+urv-wq \overline{w} \\
\end{pmatrix},
\\ u,v,w,\overline{w} \in \mathbf{Z}, \ \ u \overline{w}-vw=\pm1 \tag{2.12} 
\end{gathered} 
\end{equation}
leads to an equivalent regular $\Gamma(N,p',q')$ lattice, and we get finitely many affine lattice equivalence classes (whose number
depends on the fixed $N$ by the 

\begin{equation}
\begin{gathered}
\text{divisors of} \ \ p(N-p)-1 \ \ \text{for possible} \ q \ \text{and} \ r \ \text{parameters}, 
\end{gathered} \tag{2.12'}
\end{equation}
as we shall see in the proof of Theorem 5 and in Theorem 6 (for $N=6$). $\square$

Thus we obtain by the above discussions the following summary:  
\begin{theorem} Translations $\tau_1, \tau_2, \tau_3$ generate a lattice $\Gamma(\Phi)(\tau_1,\tau_2,\tau_3)$  
in $\SOL$ space 
if and only if $\tau_1, \tau_2$ in 
the $(O,x,y)$ plane generate a lattice $\Gamma^0$, affinely so unimodularly equivalent to a $\Gamma(N,p,q)$ lattice; 
Moreover, $\SOL$ lattices form two main classes I,~II by Theorem 3, 

In the above cases the parameter $t_3^3$ 
in $\tau_3$ satisfies the following equations and equivalence:
\begin{equation}
     \begin{gathered}
     2 \cosh{t_3^3}=p+s=N, \\
     \begin{pmatrix}
         e^{-t_3^3}&0 \\
         0&e^{t_3^3}
     \end{pmatrix}\Leftrightarrow \begin{pmatrix}
              \cosh{t_3^3}&\sinh{t_3^3} \\
              \sinh{t_3^3}&\cosh{t_3^3}
       \end{pmatrix}, \ i.e. \ t_3^3=\log(\frac{1}{2}({N+\sqrt{N^2-4})}). \square  
       \end{gathered} \tag{2.13}
     \end{equation}
\end{theorem}
By Definition 3 and by the above considerations, the affine equivalence classes of $\SOL$ lattices can be defined.
\begin{definition}
A lattice $\Gamma^1(\Phi^1)$ in $\SOL$ space is called affinelly equivalent to a lattice $\Gamma^2(\Phi^2)$ $(\Gamma^1(\Phi^1)$ $\sim \Gamma^2(\Phi^2))$, 
if there is an affine transformation $\alpha$ with
$\Gamma^2(\Phi^2)=\Gamma^1(\Phi^1) \alpha$, so that any $\SOL~$-~symmetry ($\SOL$-~isometry) $S^\alpha$ of $\Gamma^2(\Phi^2)$ is derived by a symmetry $S$ of $\Gamma^1(\Phi)$ 
with conjugacy $S^ \alpha = \alpha^{-1} S \alpha$. $\square$
\end{definition}

\begin{remark} This definition admits infinite equivalence classes, depending on $N>2$ at least, which can be expressed again by specified unimodular conjugacy, analogously as above.
In Section 3 we shall introduce for a fixed $N$ 
another algorithmic classification which leads to the so-called Bravais lattice types in $\SOL$ space. Each type, for a given $N$, shall have finitely many
affine classes, according to (2.12), (2.12') left for an "elementary" geometric number theoretical study on $N$ in the future.
\end{remark}
\subsection{$\SOL$ "parallelepipeds"}
If we take integers as coefficients, then we generate the discrete group $ \langle \tau_1,\tau_2, \tau_3 \rangle $ 
denoted by $\Gamma(\Phi)$, as above in Theorem 2. 

We know that $\SOL /\Gamma(\Phi)$ is a compact manifold, i.e. a $\SOL$ space form.

Let $\widetilde{\mathcal{F}}$ be a {\it fundamental domain} of $\Gamma(\Phi)$ with face identifications on its pairs of bent side faces (Fig.~2). 
The homogeneous coordinates of vertices of $\widetilde{\mathcal{F}}$ can be determined in our affine model by the translations
in Definition 1 with the parameters $t_i^j, \ i\in\{1,2,3\}, \ j \in \{1,2,3\}$ as follows for the cases I and II in Theorem 3 (see Fig.~2). The skew curves
$PP^{\tau_3},P'P'^{\tau_3}, QQ^{tau_3}$ are so-called translation curves (e.g. in \cite{BSz},\cite{MoSz}), defined by a first order differential equation,
expressing that any tangent of the curve can be obtained as translate of the initial tangent. Thus the bent side faces of $\widetilde{\mathcal{F}}$ can 
explicitely computed and visualized.
\begin{equation}
\begin{gathered} 
P(1,t_1^1,t_1^2,0), \ P'(1,t_2^1,t_2^2,0), \ P_3(1,t_3^1,t_3^2, t_3^3), \ Q(1,t_1^1+t_2^1,t_1^2+t_2^2,0), \\ 
Q'(1,(t_1^1+t_2^1) e^{-t_3^3}, (t_1^2+t_2^2) e^{t_3^3} ,0), \\
Q^{\tau_3}(1,t_3^1+(t_1^1+t_2^1) e^{-t_3^3}, t_3^2+(t_1^2+t_2^2) e^{t_3^3} ,t_3^3), \
P''(1,t_2^1 e^{-t_3^3}, t_2^2 e^{t_3^3} ,0), \\ P'^{\tau_3}(1,t_3^1+t_2^1 e^{-t_3^3}, t_3^2+t_2^2 e^{t_3^3} ,t_3^3), P^{\tau_3}(1,t_3^1+t_1^1 e^{-t_3^3}, 
t_3^2+t_1^2 e^{t_3^3} ,t_3^3). \tag{2.14}
\end{gathered}
\end{equation}

In Fig.~2 the case I is illustrated for parameters $N=p+s=3$ and $t_1^1=\frac{1}{\sqrt{2}}=t_1^2$, $t_2^1=\frac{3-\sqrt{5}}{2\sqrt{2}}, t_2^2=\frac{3+\sqrt{5}}{2\sqrt{2}}$,
$t_3^1=t_3^2=0, \ t_3^3=\log{\frac{3+\sqrt{5}}{2}}$.
\begin{figure*}[]
\includegraphics[width=0.9\textwidth]{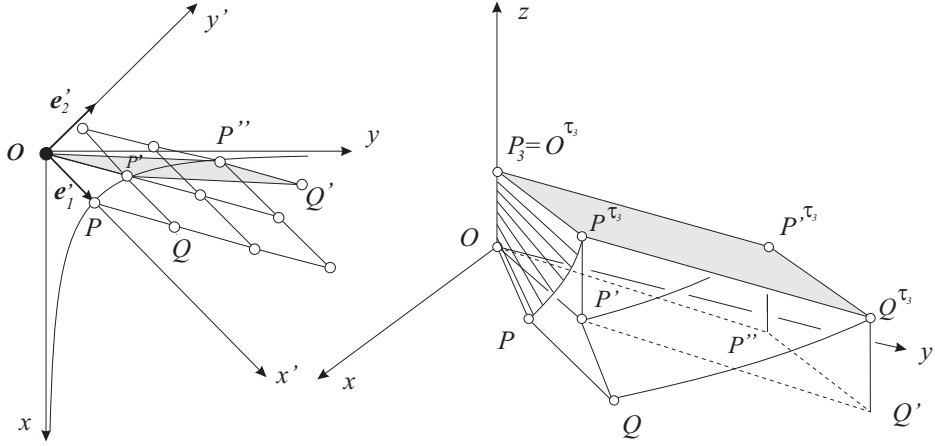} 
\caption{A fundamental parallelepiped for a "fundamental lattice" $N=3$, $p=0$, $s=3$, $q=1$, $r=-1$.}
\label{Fig:2}
\end{figure*}
\begin{definition}
A $\SOL$ point lattice $\Gamma_P(\Phi)$ is a discrete orbit of point $P$ in $\SOL$ space, generated by an arbitrary lattice $\Gamma(\Phi)$ above.
For visualizing a point lattice we have chosen the origin as starting point, by homogeneity of $\SOL$. \ $\square$
\end{definition}
\section{Bravais lattice types in $\SOL$ space}
\label{sec:3}
In Euclidean geometry and crystallography a Bravais point lattice, named to honor of Auguste {\sc Bravais} (1811-1863), 
is an infinite set of points generated by discrete "translation operations". The number of Bravais lattices (or lattice types)
in three-dimensional Euclidean space is well known to be 14 if, as usual, a lattice type is defined by its lattice isometry group up to affine
(so unimodular by basis change) conjugacy. 

Similarly to Euclidean geometry we introduce the notion of Bravais lattice type in $\SOL$ space, but now this will be coarser than affine (unimodular)
conjugacy.

\begin{definition}
Let the group of all isometries in $\SOL$ space, fixing the origin, be denoted by $O_{\SOL}$ (see (1.5), (1.6)). 
Its subgroup is called a {\it point group} of $\SOL$ space, if the group leaves a lattice invariant. We consider an arbitrary point lattice $\Gamma_O(\Phi)$,
its (maximal) point group is denoted by
$S_O(\Gamma(\Phi))$. \ $\square$
\end{definition}
\begin{definition}[Bravais lattice type]
$\SOL$ point lattices $\Gamma_O(\Phi)=\langle \tau_1,\tau_2,\tau_3 \rangle$ by fixed $N>2$ above will be {\it classified in an algorithmic way, as concretely follows}
at Theorem 5: {\it first}, by Theorem 3 into main cases I and II; {\it second}, by its maximal point groups $S_O(\Gamma(\Phi))$;
{\it third}, by the existence of a linear $z$-sublattice $\Gamma_{Oz} \subset \Gamma_O(\Phi)$; {\it fourth}, by the "centering" of the {\it base-lattice}
$\Gamma^0(\Phi)$ and {\it its parallel} $\Gamma^{t_3^3}$, induced by $\tau_3(t_3^1,t_3^2,t_3^3)$. \ $\square$ 
\end{definition}
That means, affine equivalent $\SOL$ lattices will be in the same Bravais type, but not vice-versa, in general (see in Theorem 6). Moreover, we allow that, depending on $N$,
 some Bravais types do not appear for $N$. That means, we also have essential differences in comparison with Euclidean geometry.
\begin{theorem}
The lattices, in main cases I and II by Theorem 3 for fixed $N$, can form 17 (seventeen) different Bravais types in $\SOL$ space.
\end{theorem}
{\it Proof:}

We classify the lattices, by their point groups in the following.
The possible subgroups, of $D_4 \sim O_{\SOL}$ will be listed first by (1.5), we see conjugate subgroups by $O_{\SOL}$ besides each other.
\begin{enumerate}
\item Group of order 1: $\mathbf{Id}$,
\item Groups of order 2:\\
$$\mathbf{D_x}: \Bigg\{\mathbf{Id}, ~ \begin{pmatrix}
         1 & 0 & 0 \\
         0 & -1 & 0 \\
         0 & 0 & 1
     \end{pmatrix} \Bigg\}, \ 
\mathbf{D_y}: \Bigg\{  \mathbf{Id}, ~ \begin{pmatrix}
         -1 & 0 & 0 \\
         0 & 1 & 0 \\
         0 & 0 & 1
     \end{pmatrix} \Bigg\},$$ \\
$$\mathbf{C_2}: \Bigg\{\mathbf{Id}, ~ \gamma_2=\gamma_2^{-1}=\begin{pmatrix}
              -1 & 0 & 0 \\
              0 & -1 & 0 \\
              0 & 0 & 1
     \end{pmatrix} \Bigg\}, $$ \\ 
$$\mathbf{D_r}: \Bigg\{ \mathbf{Id}, ~ \delta_r=\delta_{r}^{-1}=\begin{pmatrix}
         0 & 1 & 0 \\
         1 & 0 & 0 \\
         0 & 0 & -1
     \end{pmatrix} \Bigg\},
\mathbf{D_{\overline{\mathbf{r}}}}: \Bigg\{ \mathbf{Id}, ~ \delta_{\overline{r}}=\delta_{\overline{r}}^{-1}=\begin{pmatrix}
         0 & -1 & 0 \\
         -1 & 0 & 0 \\
         0 & 0 & -1
     \end{pmatrix} \Bigg\}.$$     
\item Groups of order 4:\\
$$\mathbf{D_2}: \Bigg\{\mathbf{Id},~ 
\begin{pmatrix}
              -1 & 0 & 0 \\
              0 & -1 & 0 \\
              0 & 0 & 1
     \end{pmatrix}; 
     \begin{pmatrix}
         -1 & 0 & 0 \\
         0 & 1 & 0 \\
         0 & 0 & 1
     \end{pmatrix}, \begin{pmatrix}
         1 & 0 & 0 \\
         0 & -1 & 0 \\
         0 & 0 & 1
     \end{pmatrix}   \Bigg\},$$
$$\mathbf{D_{\overline{2}}}: \Bigg\{\mathbf{Id},~  \gamma_2=\begin{pmatrix}
         -1 & 0 & 0 \\
         0 & -1 & 0 \\
         0 & 0 & 1
     \end{pmatrix}, ~ $$
     $$
     \delta_r=\begin{pmatrix}
              0 & 1 & 0 \\
              1 & 0 & 0 \\
              0 & 0 & -1
     \end{pmatrix},~
     \delta_{\overline{r}}=\begin{pmatrix}
              0 & -1 & 0 \\
              -1 & 0 & 0 \\
              0 & 0 & -1
     \end{pmatrix}=\gamma_2 \delta_r=\delta_r \gamma_2 \Bigg\},$$
$$\mathbf{C_4}: \Bigg\{\mathbf{Id},~ \gamma_4=\begin{pmatrix}
         0 & 1 & 0 \\
         -1 & 0 & 0 \\
         0 & 0 & -1
     \end{pmatrix}, ~ \gamma_2=\begin{pmatrix}
         -1 & 0 & 0 \\
         0 & -1 & 0 \\
         0 & 0 & 1
     \end{pmatrix}, ~ \gamma_4^{-1}=\begin{pmatrix}
              0 & -1 & 0 \\
              1 & 0 & 0 \\
              0 & 0 & -1
     \end{pmatrix}\Bigg\}.$$ \\
\item The whole group of order 8:\\
$$\mathbf{D_4}: \Bigg\{\mathbf{Id},~ 
     \begin{pmatrix}
              0 & 1 & 0 \\
              -1 & 0 & 0 \\
              0 & 0 & -1
     \end{pmatrix},~
     \begin{pmatrix}
             -1 & 0 & 0 \\
              0 & -1 & 0 \\
              0 & 0 & 1
     \end{pmatrix},~ 
     \begin{pmatrix}
         0 & -1 & 0 \\
        1 & 0 & 0 \\
         0 & 0 & -1
     \end{pmatrix}, 
      $$ 
$$ \begin{pmatrix}
         -1 & 0 & 0 \\
         0 & 1 & 0 \\
         0 & 0 & 1
     \end{pmatrix}, \begin{pmatrix}
         1 & 0 & 0 \\
         0 & -1 & 0 \\
         0 & 0 & 1
     \end{pmatrix},
     \begin{pmatrix}
         0 & 1 & 0 \\
         1 & 0 & 0 \\
         0 & 0 & -1
     \end{pmatrix}, \begin{pmatrix}
         0 & -1 & 0 \\
        -1 & 0 & 0 \\
         0 & 0 & -1
     \end{pmatrix}\Bigg\}.$$
\end{enumerate}

By Theorem 3 $\SOL$ lattices form two main classes which will be refined by Definition 8 yielding the Bravais lattice types.

We have by \cite{AQ} the following 
\begin{lemma} 
Base lattice $\Gamma^0(\Phi)$ of a $\SOL$-lattice $\Gamma(\Phi)$ 
does not contain a non-zero isotropic (light-like) translation vector proportional either to $\be_1$ or to $\be_2$. 
\end{lemma}
{\it Proof:}

The symmetry group $S(\Gamma(\Phi))$ of $\Gamma(\Phi)$ contains a hyperbolic rotation $H(O,t_3^3)$ with centre $O$ and by parameter $t_3^3$. 
If $\mathbf{v} \in \Gamma^0(\Phi)$ an isotropic vector then its images would accumulate either at $H^n$ or at $H^{-n}, \ (n\in \mathbf {N})$, since
\[\lim_{n \to \infty}{e^{-nx} \rightarrow 0} \ \ (\text{for} \ 0<x\in \mathbf{R}) \] in contradiction
to discreteness. $\square$
\begin{corollary}
In main cases I and II the point group of a lattice does not admit reflections in the $(y,z)$ and $(x,z)$ planes, respectively:
$$\mathbf{D_y}=\begin{pmatrix}
         -1 & 0 & 0 \\
         0 & 1 & 0 \\
         0 & 0 & 1
     \end{pmatrix},\ \ \mathbf{D_x}=\begin{pmatrix}
         1 & 0 & 0 \\
         0 & -1 & 0 \\
         0 & 0 & 1
     \end{pmatrix}. \square$$
\end{corollary}
Else the lattice would contain isotropic translation in contradiction to Lemma 1.
\begin{lemma} 
{\bf In case I:} The maximal point group of every lattice has the halfturn $\gamma_2$ about axis $z$. Thus we have the possibilities: 
\begin{equation}
\mathbf{C_2},~ \mathbf{D_{\overline{2}}},~ \mathbf{C_4} \tag{3.1}
\end{equation}
that map $\tau_3(0,0,t_3^3)$ and its $z$-sublattice $\Gamma_{Oz}$ into itself or into its inverse and the base lattice $\Gamma^0(\Phi)$ will be invariant (Fig.~3-5).

{\bf In case II:} The possible maximal point groups are
\begin{equation}
\mathbf{Id}, ~ \mathbf{C_2},~ \mathbf{D_r} (\sim \mathbf{D_{\overline{r}}}),~ \mathbf{D_{\overline{2}}},~ \mathbf{C_4}, \tag{3.2}
\end{equation}
leaving $\Gamma^0(\Phi)$ and $\Gamma(\Phi)$ invariant.
\end{lemma}

{\it Proof for case I:}

The symmetry of lattice $\Gamma(\Phi)$ under $\mathbf{C_2}$ is clear in this case.

Conjugacy by $\delta_r \in \mathbf{D_r} \subset \mathbf{D_{\overline{2}}}$ produces $\tau_2=\tau_1^{\delta_r}=\delta_r\tau_1\delta_r$ $\sim 
(t_2^1=t_1^2,~t_2^2=t_1^1,0)$ as possible second basis translation, if $t_1^2 \ne t_1^1$ and $t_1^2 \ne -t_1^1$. Then translation 
$\tau_1\tau_1^{\delta_r}
=:\tau_1^,$ points into direction of symmetry axis $r$ of axes $x$ and $y$; $\tau_2^,:=\tau_1^{-1}\tau_1^{\delta_r} \sim (-t_1^1+t_1^2, -t_1^2+t_1^1,0)$ 
points into direction of second symmetry axis $\overline{r}$ of $x$ and $y$ (Fig.~4.a). Besides {\it the first, centred base lattice} before, 
$\tau_1^,(t_1^1,t_1^1)$ and $\tau_2^,(-t_2^2,t_2^2)$
can be basis of {\it the second, primitive lattice} $\Gamma^0(\Phi)$ (Fig~.~4.a-b, like in Euclidean plane $\mathbf{E}^2$). We always exclude isotropic translation vectors in the following.

At $\mathbf{C_4}$: basis translation $\tau_1(t_1^1,t_1^2,0)$ involves its $\gamma_4$-conjugate $\tau_1^{\gamma_4}=\gamma_4^{-1} \tau_1 \gamma_4 \sim
(-t_1^2,$ $t_1^1,0)$ as possible second basis translation $\tau_2^,$ in $\Gamma^0(\Phi)$ (Fig.~5), if $t_1^2/t_1^1 \notin {\mathbf{Q}}$ (rational numbers). Else isotropic
translation would occur.

{\it Proof for case II:}

At $\mathbf{C_2}$: basis translation $\tau_3(t_3^1,t_3^2,t_3^3)$ implies that 
\begin{equation}
\begin{gathered}
(\tau_3^{\gamma_2})=\gamma_2 \tau_3 \gamma_2 \sim (-t_3^1,-t_3^2,t_3^3), \\
(\tau_3^{\gamma_2})^{-1}=\gamma_2 \tau_3^{-1}\gamma_2 \sim (t_3^1 e^{t_3^3}, t_3^2 e^{-t_3^3},-t_3^3) \ \ \text{and e.g.} \\
(\tau_3^{\gamma_2})^{-1}\tau_3=\gamma_2 \tau_3^{-1}\gamma_2 \tau_3 \sim (2t_3^1,2t_3^2,0) \in \Gamma^0(\Phi), \\
(\tau_3^{\gamma_2})\tau_3=\gamma_2\tau_3\gamma_2\tau_3 \sim (t_3^1(1-e^{-t_3^3}), \  t_3^2(1-e^{t_3^3}), \  2t_3^3) \nparallel (0,0,z),
\end{gathered} \tag{3.3}
\end{equation}
i.e. there is no sublattice in $z$-direction in this case. For $n\in \bN$ we have
\begin{equation}
\begin{gathered}
\tau_3^n \sim \big( t_3^1 \Big(\frac{1-e^{-n t_3^3}}{1-e^{-t_3^3}}\Big), \  t_3^2 \Big(\frac{e^{n t_3^3}-1}{e^{t_3^3}-1}\Big), \ nt_3^3\big), \\
\tau_3^{-n} \sim \big( -t_3^1 e^{t_3^3} \Big(\frac{e^{(n-1) t_3^3}-1}{e^{t_3^3}-1}\Big), \  
-t_3^2 e^{-t_3^3} \Big(\frac{e^{-(n-1) t_3^3}-1}{e^{-t_3^3}-1}\Big), \ -nt_3^3\big).
\end{gathered} \tag{3.4}
\end{equation}
This means in homogeneous coordinates that $(0,0,1,0)$, the ideal point of $y$-axis, is the limit if $n \rightarrow \infty$; and $(0,1,0,0)$,
the ideal point of $x$-axis, is the limit if $n \rightarrow -\infty$.

At $\mathbf{D_r}$: for $\tau_3$ we have by $\delta_r$-conjugacy
\begin{equation}
\begin{gathered}
\tau_3^{\delta_r}=\delta_r \tau_3\delta_r \sim (t_3^2, t_3^1,-t_3^3), \
(\tau_3^{\delta_r})^{-1} \sim (-t_3^2 e^{-t_3^3}, -t_3^1 e^{t_3^3},t_3^3) \ \ \text{and e.g.} \\
(\tau_3^{\delta_r})\tau_3=\delta_r \tau_3\delta_r \tau_3 \sim (t_3^2 e^{-t_3^3}+t_3^1, t_3^1 e^{t_3^3}+t_3^2,0) \in \Gamma^0(\Phi), \\
[(\tau_3^{\delta_r})\tau_3]^{\delta_r}=\tau_3\delta_r \tau_3\delta_r  \sim (t_3^1 e^{t_3^3}+t_3^2, t_3^2 e^{-t_3^3}+t_3^1,0) \in \Gamma^0(\Phi), \\ 
(\tau_3^{\delta_r})^{-1}\tau_3=\delta_r \tau_3^{-1}\delta_r \tau_3 \sim (-t_3^2 e^{-2t_3^3}+t_3^1, -t_3^1 e^{2t_3^3}+t_3^2,2t_3^3).
\end{gathered} \tag{3.5}
\end{equation}
That means, $(\tau_3^{\delta_r})\tau_3$ points neither in $r$-dirction, nor in $\overline{r}$-direction, and it can be zero iff $t_3^2=-t_3^1 e^{t_3^3}$.
Moreover, $(\tau_3^{\delta_r})^{-1}\tau_3$ points in $z$-direction iff $t_3^2=t_3^1 e^{2t_3^3}$.

At $\mathbf{D_{\overline{2}}}$: besides $\delta_r$-conjugacy as at $\mathbf{D_r}$ in formulas (3.5), here analogous $\delta_{\overline{r}}$-conjugacy,
or rather preferred, $\gamma_2$-conjugacy at $\mathbf{C_2}$ 
also occurs:
\begin{equation}
\begin{gathered}
(\tau_3^{\gamma_{2}})^{-1} \tau_3 \sim (2t_3^1, 2t_3^2,0) \in \Gamma^0(\Phi). 
\end{gathered} \tag{3.6}
\end{equation}

At $\mathbf{C_4}$: for $\tau_3$ we have by $\gamma_4$-conjugacy
\begin{equation}
\begin{gathered}
\tau_3^{\gamma_{4}} \sim (-t_3^2, t_3^1,-t_3^3), \ \ 
\tau_3^{\gamma_{4}^{-1}} \sim (t_3^2, -t_3^1, -t_3^3), \\
(\tau_3^{\gamma_{4}})^{-1} \sim (t_3^2 e^{-t_3^3}, -t_3^1 e^{t_3^3},t_3^3), \ \
(\tau_3^{\gamma_{4}^{-1}})^{-1} \sim (-t_3^2 e^{-t_3^3}, t_3^1 e^{t_3^3},t_3^3), \ \ \text{and e.g.} \\
(\tau_3^{\gamma_{4}})\tau_3=\gamma_{4}^{-1} \tau_3\gamma_{4} \tau_3 \sim (-t_3^2 e^{-t_3^3}+t_3^1, t_3^1 e^{t_3^3}+t_3^2,0) \in \Gamma^0(\Phi), \\
(\tau_3^{\gamma_{4}})^{-1}\tau_3=\gamma_{4}^{-1} \tau_3^{-1}\gamma_{4} \tau_3 \sim (t_3^2 e^{-2t_3^3}+t_3^1, -t_3^1 e^{2t_3^3}+t_3^2,2t_3^3), \\
(\tau_3^{\gamma_{2}})^{-1}\tau_3 \sim (2t_3^1, 2t_3^2,0),
\end{gathered} \tag{3.7} 
\end{equation} 
That means $(\tau_3^{\gamma_{4}})\tau_3$ is isotropic if $t_3^2=t_3^1 e^{t_3^3}$ or $t_3^2=-t_3^1 e^{t_3^3}$, these will be excluded.
We have no translation into $z$-direction. Furthermore, $\gamma_2$-conjugacy of $\tau_3$ also occurs, as formerly.
Thus we have proved Lemma 2 with observations on lattice $\Gamma(\Phi)$, will be used later on. $\square$
\medbreak
For visualization only, we have indicated the Euclidean size of angles between the basis vectors $\tau_1,~\tau_2,~\tau_3$ by 
$\alpha=(\tau_2,\tau_3)\angle,~\beta=(\tau_3,\tau_1)\angle$, $\gamma=(\tau_1,\tau_2)\angle$ in Fig.~3-9. 
In analogy to Euclidean geometry, Bravais types of $\SOL$-lattices can be discussed and we have the following cases:
\medbreak
{\bf{Lattice types in our main case I: $t_1^3=t_2^3=t_3^1=t_3^2=0$}}
\medbreak
\begin{enumerate}
\item[$\circ$]  {\it The point group is $\mathbf{C_2}$}.

{\bf I/1}: Now, we may consider that case where , $\alpha=\beta=\frac{\pi}{2}$, $0<\gamma<\frac{\pi}{2}$. 
This is the type of primitive monoclinic lattice (see Fig.~3, in analogy to $\EUC$).
\begin{figure*}[ht]
\includegraphics[width=10cm]{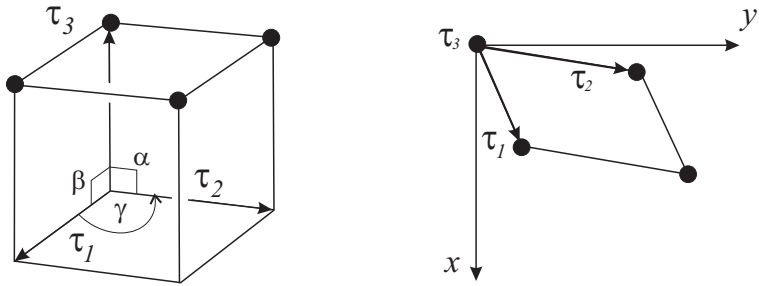} 
\caption{$\mathbf{C_2}$, \ {\bf I/1} Primitive monoclinic lattice}
\label{Fig:3}
\end{figure*}
\item[$\circ$]  {\it The point group is $\mathbf{D_{\overline{2}}}$}.

{\bf I/2}: Now, we consider that case, as before but $\tau_1^{\delta_r}=\tau_2$ is chosen, i.e. 
\begin{equation}
\begin{gathered}
\begin{pmatrix}
         t_1^1 & t_1^2 \\
         t_2^1 & t_2^2 
\end{pmatrix}= \begin{pmatrix}
         t_1^1 & t_1^2 \\
         t_1^2 & t_1^1 
\end{pmatrix}, \ \ \text{and e.g.} \ t_1^1>t_1^2>0 \ \ N \ge 3, \ \ p=0, \ s=N, \\ q=1, \ r=-1 \ \ \text{by (2.9)}. \notag 
\end{gathered}
\end{equation}
Translations in directions to symmetry axes $r$ and 
$\overline{r}$ of $x$- and $y$-axes are guaranteed (see Fig.~4.a, like a base centred orthorhombic lattice in $\mathbf{E}^3$). 
\begin{figure*}[ht]
\includegraphics[width=10cm]{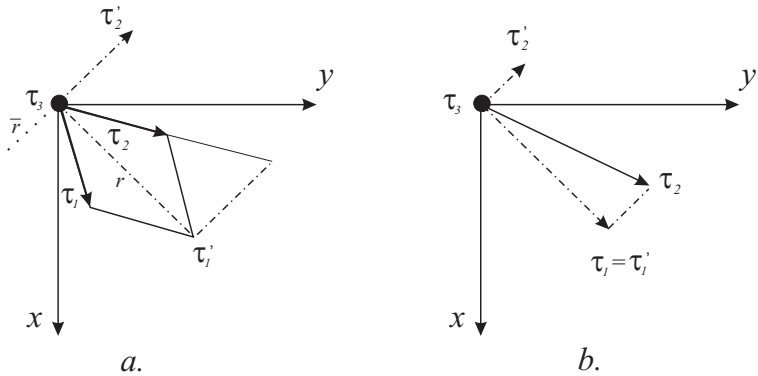} 
\caption{$\mathbf{D_{\overline{2}}}$, \ \  a. {\bf I/2} Base centred orthorhombic lattice; \ \ 
b. {\bf I/3} Primitive orthorhombic lattice}
\label{Fig:4}
\end{figure*}
\medbreak
{\bf I/3}: Now, we consider that case where $\alpha=\beta=\gamma=\frac{\pi}{2}$ and basis vectors $\tau_1^,(t_1^1>0,t_1^2=t_1^1,0), \ \tau_2^,(t_2^1=-t_2^2,t_2^2>0,0)$ 
are in directions of symmetry axes of $x$- and $y$-axis, respectively. We assume $t_2^2/t_1^1 \notin \mathbf{Q}$,
else isotropic translation would occur. But a new basis vector $\tau_2$ can also be chosen in positive $(x,y)$ domain.
That means, we have obtained a {\it primitive orthorhombic lattice} (Fig.~4.b, as an analogue to that in $\mathbf{E}^3$).
\item[$\circ$]  The point group is $\mathbf{C_4}$.

{\bf I/4}: Now, we consider that case where $\alpha=\beta=\gamma = \frac{\pi}{2}$, and 
$\tau_{2}^,= \tau_1^{\gamma_4}:=\gamma_4^{-1} \tau_1 \gamma_4 \sim (t_2^1=-t_1^2, \ t_2^2=t_1^1,0)$ is obtained by $\gamma_4$ conjugacy from 
$\tau_1 \sim (t_1^1>0,t_1^2>0,0)$ (and e.g. 
$\tau_2=\tau_1+\tau_{2}^,=\tau_1\tau_{2}^,=\tau_1\gamma_4^{-1}\tau_1\gamma_4$ and $\tau_1$ both point into the positive domain of $(x,y)$). 
This is a {\it primitive "tetragonal"} lattice in Euclidean
sense (Fig.~5). We exclude isotropic translation, i.e. $t_1^1/t_1^2 \notin \mathbf{Q}$ (rationals). Furthermore, (2.9) implies now
\begin{equation}
\frac{N-2p+\sqrt{N^2-4}}{2q}=\frac{2q}{\sqrt{N^2-4}-N+2p}, \ \ \text{i.e.} \ \ qq=p(N-p)-1 \tag{3.7'}
\end{equation}
for $\mathbf{C_4}$, as indicated at (2.12')

In main case I we have easily obtained 4 Bravais types {\bf I/1-I/4}.
\begin{figure*}[ht]
\includegraphics[width=10cm]{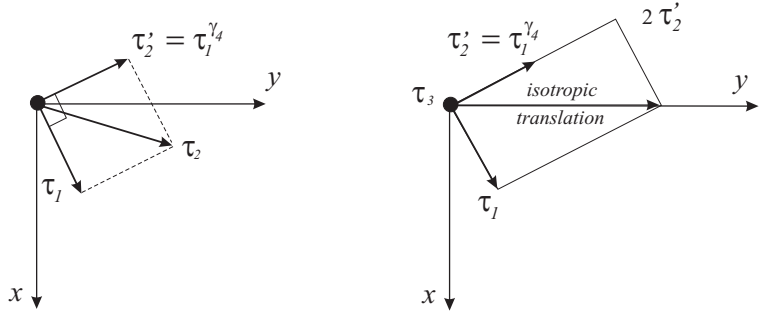} 
\caption{{\bf I/4}, \ $\mathbf{C_4}$ Primitive tetragonal lattice (without isotropic translation)}
\label{Fig:5}
\end{figure*}
\end{enumerate}

{\bf{Lattice types in our main case II:}} $t_1^3=t_2^3=0, \ (0,0) \ne (t_3^1,t_3^2)$ can be chosen $\text{mod}(\Gamma^0 \langle \tau_1,\tau_2 \rangle)$.

\begin{enumerate}
\item[$\circ$]  {\it The point group is trivial}. 

{\bf II/1}: This type has one affine class for fixed $N$ the so-called {\it primitive triclinic lattice}. 
(see Fig.~6, like in $\mathbf{E}^3)$.
\begin{figure*}[ht]
\includegraphics[width=10cm]{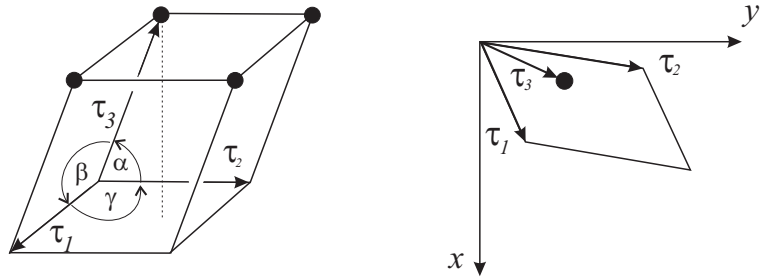} 
\caption{{\bf Id}, \ {\bf II/1.} Primitive triclinic lattice}
\label{Fig:6}
\end{figure*}
\item[$\circ$]  {\it The point group is $\mathbf{C_2}$}.

{\bf II/2}: In this case, by Lemma 2, $(\tau_3^{\gamma_2})^{-1}\tau_3 \in \Gamma^0(\Phi)$, i.e.  $\gamma_2{\tau_3}^{-1}\gamma_2 \tau_3=2{\tau_3}^,=
\tau_1$ or $\tau_2$ or $\tau_1+\tau_2 $) can be chosen for the projected translation of $\tau_3$ 
into $[x,y]$ plane) (Fig.~7, formulas (3.3)). We have got an analogue of {\it monoclinic single-side-face-centred lattice}, i.e. 
we have one additional lattice point pair at the centres of 
some opposite side faces; or in a new parallelepiped (with new basis) if $2 \tau_3^,=\tau_1+\tau_2 $ by usual convention (Fig.~7.a-b).
All these subcases are affinely equivalent, but not equivalent with the "primitive" case {\bf I/1} where a sublattice $\Gamma_{Oz}$ 
exists on $z$-axis. Here
such a sublattice does not exist, since neither $(\tau_3^{\gamma_2})\tau_3 \sim (t_3^1(1-e^{-t_3^3}), t_3^2(1- e^{t_3^3}), 2 t_3^3)$ nor $\tau_3\tau_3$,
etc. have zero first and second components at the same time. 

{\it Thus the analogy to Euclidean centering is not exact in general}.
\medbreak
\begin{figure*}[ht]
\includegraphics[width=10cm]{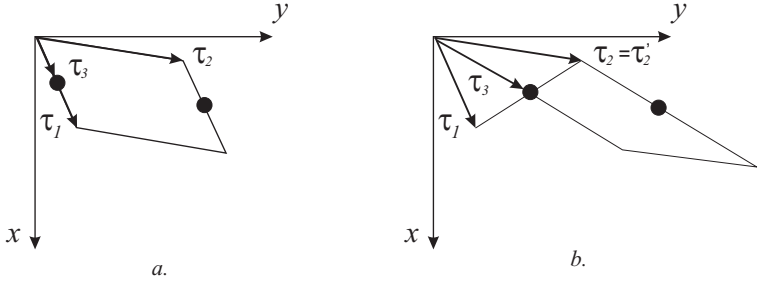} 
\caption{$\mathbf{C_2}$, \ {\bf II/2.} Monoclinic single-side-face-centred lattice}
\label{Fig:7}
\end{figure*}
\item[$\circ$]  {\it The point group is $\mathbf{D_r}$} (Fig.~8, formulas (3.5), see also details at $\mathbf{D_{\overline{2}}}$ as follows).
\begin{enumerate}
\item[i.]Assume that $t_3^2/t_3^1=:\overline{x}= e^{2t_3^3}$, i.e. $(\tau_3^{\delta_r})^{-1}\tau_3 \sim (0,0,2t_3^3)$ points to $z$-direction (similarly to 
$t_3^2/t_3^1=:\overline{x}= -e^{2t_3^3}$ for $\delta_{\overline{r}}$). 
\\ 
Then $(\tau_3^{\delta_r})\tau_3 \sim (t_3^1(e^{t_3^3}+1), t_3^1 e^{t_3^3}(1+e^{t_3^3}),0) \in \Gamma^0(\Phi)$, not on $r$ and $\overline{r}$
(Fig.~8.a).
\item[i.a.] {\bf II/3}. Assume that $t_3^2=t_3^1 e^{2t_3^3}$ and the base lattice $\Gamma^0(\Phi)=\langle \tau_1, \tau_2 \rangle$ is centred, and 
say $(\tau_3^{\delta_r})\tau_3=\tau_2=\tau_1^{\delta_r},
 \ \text{i.e.} \ \tau_1 \sim (t_1^1=t_3^1 e^{t_3^3} (e^{t_3^3}+1), t_1^2=t_3^1(e^{t_3^3}+1),0)$ (Fig.~8.a). 
 Moreover, we know that $e^{t_3^3} \notin \mathbf{Q}$, 
 to exclude isotropic translation. 
\item[i.b.] {\bf II/4}. Assume as above that $t_3^2=t_3^1 e^{2t_3^3}$ and the base lattice $\Gamma^0(\Phi)=\langle \tau_1, \tau_2 \rangle$ is
"primitive", i.e. $(\tau_3^{\delta_r})\tau_3=\tau_1^,\tau_2^,$ where
$\tau_1^,$ is in direction $r, \tau_2^,$ is on $\overline{r}$ (Fig.~8.b). That means, $\tau_1^, \sim (t_1^1=t_3^1 (e^{t_3^3}+1)^2, t_1^2=t_1^1=t_3^1(e^{t_3^3}+1)^2,0)$,
$\tau_2^, \sim (t_2^1=-t_2^2=-t_3^1 (e^{2t_3^3}-1), t_2^2=t_3^1(e^{2t_3^3}-1),0)$. Again, we exclude isotropic translation.
\item[ii.c.] {\bf II/5}. If $\pm t_3^1 e^{2t_3^3}\ne t_3^2 \ne \mp t_3^1 e^{2t_3^3}$, then $(\tau_3^{\delta_r})^{-1}\tau_3$ does not point into direction of $z$-axis, but
$(\tau_3^{\delta_r})\tau_3 \sim (t_3^2 e^{-t_3^3}+t_3^1, t_3^1 e^{t_3^3} +  t_3^2,0) \in \Gamma^0(\Phi)=\langle \tau_1, \tau_2 \rangle$ as before.
First, say $(\tau_3^{\delta_r})\tau_3=\tau_2=\tau_1^{\delta_r}$, i.e. the base lattice $\Gamma^(0)(\Phi)$ is centred, without isotropic translation (Fig.~8.c).
\item[ii.d.] {\bf II/6}. Then comes $\pm t_3^1 e^{2t_3^3}\ne t_3^2 \ne \mp t_3^1 e^{2t_3^3}$, but $(\tau_3^{\delta_r})\tau_3=\tau_1^,\tau_2^,$, i.e. $\Gamma^0(\Phi)$ is primitive, without isotropic translation
(Fig.~8.d).
\item[iii.e.] {\bf II/7}. Assume $t_3^2=-t_3^1 e^{t_3^3}$ as special subcase, where $(\tau_3^{\delta_r})\tau_3=0$ (or $t_3^2=t_3^1 e^{t_3^3}$ 
where $(\tau_3^{\delta_{\overline{r}}})\tau_3=0$). The base lattice $\Gamma^0(\Phi)$ 
is not restricted yet by $\tau_3$. $\Gamma^0(\Phi)$ can be centred, without isotropic translation (Fig.~8.e).
\item[iii.f.] {\bf II/8.} Assume $t_3^2=-t_3^1 e^{t_3^3}$ (or $t_3^2=t_3^1 e^{t_3^3}$ 
where $(\tau_3^{\delta_{\overline{r}}})\tau_3=0$) as in previous case but with primitive base lattice $\Gamma^0(\Phi)$, without isotropic translation (Fig.~8.f).
\end{enumerate}
{\bf We remark} that at centred base lattice, formula (2.9) implies $p=0, q=1$ for given $N$. At primitive base lattice, (2.9) involve $N=2p$, i.e. $N$ is
even. That means, for odd $N$ such primitive Bravais types do not occur.
\medbreak
\item[$\circ$]  {\it The point group is $\mathbf{D_{\overline{2}}}$, i.e. $\mathbf{D_{r}}$ and $\mathbf{C_{2}}$ both together}.

Therefore we use Fig.~8 and formulas (3.5-6). Besides criteria at $\mathbf{D_{r}}$, i.e. 6 possibilities in Fig.~8, the 
additional $\mathbf{C_{2}}$ assumption (3.6) has to be fulfilled. We have in more details the following equations:
\begin{equation}
\begin{gathered}
\alpha t_1^1+\beta t_1^2= t_3^2 e^{-t_3^3}+t_3^1 \\
\alpha t_1^2+\beta t_1^1= t_3^1 e^{t_3^3}+t_3^2 \\
\gamma t_1^1+\delta t_1^2= 2 t_3^1, \ \ 
\gamma t_1^2+\delta t_1^1= 2 t_3^2 \ \ {\text{for}} \ \  
\begin{pmatrix}
t_1^1 & t_1^2 \\
t_2^1 & t_2^2 
\end{pmatrix}=
\begin{pmatrix}
t_1^1 & t_1^2 \\
t_1^2 & t_1^1 
\end{pmatrix}; \\
\text{i.e.} \ \ \alpha = \frac{t_1^1-e^{t_3^3} t_1^2}{t_1^1t_1^1-t_1^2t_1^2}(t_3^2 e^{-t_3^3}+t_3^1), \ \
\beta = \frac{t_1^1 e^{t_3^3}-t_1^2}{t_1^1t_1^1-t_1^2t_1^2}(t_3^2 e^{-t_3^3}+t_3^1); \\
\gamma= 2\frac{t_3^1 t_1^1-t_3^2 t_1^2}{t_1^1t_1^1-t_1^2t_1^2}, \ \ \delta=2\frac{t_3^2 t_1^1-t_3^1 t_1^2}{t_1^1t_1^1-t_1^2t_1^2}.
\end{gathered} \tag{3.8}
\end{equation}
hold for {\it centred base lattice} with $\alpha$, $\beta$, $\gamma$, $\delta \in \mathbf{Z}$ integers. Or similarly,
\begin{equation}
\begin{gathered}
\overline{\alpha} t_1^1-\overline{\beta} t_2^2= t_3^2 e^{-t_3^3}+t_3^1 \\
\overline{\alpha} t_1^1+\overline{\beta} t_2^2= t_3^1 e^{t_3^3}+t_3^2 \\
\overline{\gamma} t_1^1-\overline{\delta} t_2^2= 2 t_3^1, \ \ 
\overline{\gamma} t_1^1+\overline{\delta}t_2^2= 2 t_3^2 \ \ {\text{for}} \ \  
\begin{pmatrix}
t_1^1 & t_1^2 \\
t_2^1 & t_2^2 
\end{pmatrix}=
\begin{pmatrix}
t_1^1 & t_1^1 \\
-t_2^2 & t_2^2 
\end{pmatrix}, 
\end{gathered} \notag
\end{equation}
\begin{equation}
\begin{gathered}
\text{i.e.} \ \ \overline{\alpha} = \frac{1+e^{-t_3^3}}{2t_1^1}(t_3^2+ t_3^1 e^{t_3^3}), \ \
\overline{\beta} = \frac{1- e^{-t_3^3}}{2t_2^2}(t_3^2+ t_3^1 e^{t_3^3}); \\
\overline{\gamma}= \frac{t_3^1+t_3^2}{t_1^1}, \ \ \overline{\delta}=\frac{t_3^2-t_3^1}{t_2^2}.
\end{gathered} \tag{3.9}
\end{equation}
hold for {\it prmitive base lattice}, $\overline{\alpha}$, $\overline{\beta}$, $\overline{\gamma}$, $\overline{\delta} \in \mathbf{Z}$. From formulas (2.9), (2.13) come also
\begin{equation}
\begin{gathered}
\frac{t_2^1}{t_1^1}=\frac{N-2p-\sqrt{N^2-4}}{2q}, \ \ \frac{t_2^2}{t_1^2}=\frac{N-2p+\sqrt{N^2-4}}{2q},\\
e^{t_3^3}= \frac{N+\sqrt{N^2-4}}{2},\ \ e^{-t_3^3}= \frac{N-\sqrt{N^2-4}}{2}. \tag{3.10}
\end{gathered}
\end{equation}
\begin{figure*}[]
\includegraphics[width=11cm]{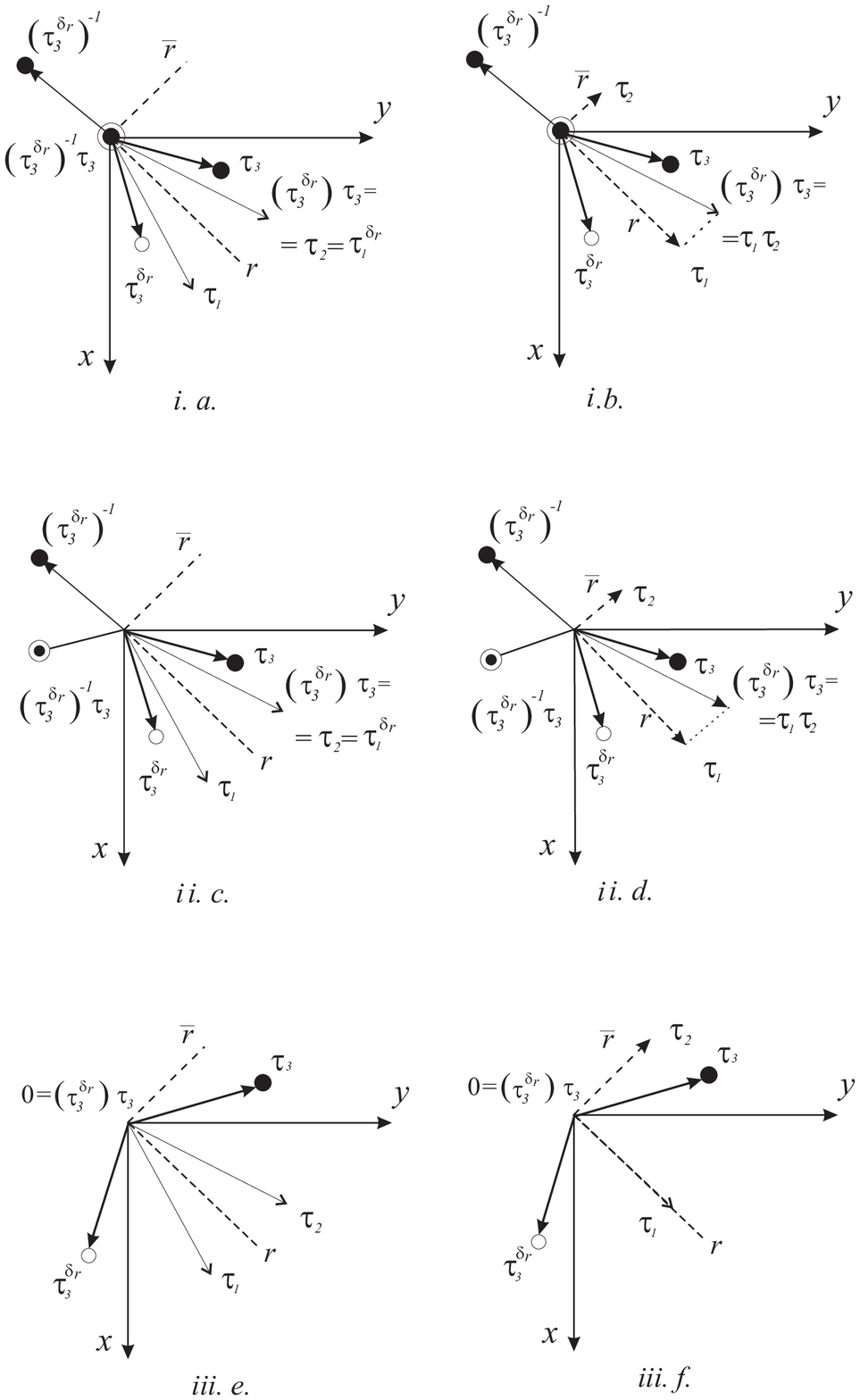} 
\caption{$\mathbf{D_r}$, Symbolic pictures, also for $\mathbf{D_{\overline{2}}}$ \newline 
i.a. {\bf II/3} Base-centred lattice with $\Gamma_{Oz}$, only with $\mathbf{D_r}$ \newline
i.b. {\bf II/4} Base-primitive lattice with $\Gamma_{Oz}$, automatically with $\mathbf{D_{\overline{2}}}$ \newline
ii.c. {\bf II/5} Base-centred lattice without $\Gamma_{Oz}$, also with $\mathbf{D_{\overline{2}}}$ for {\bf II/9} \newline
ii.d. {\bf II/6} Base-primitive lattice without $\Gamma_{Oz}$, also $\mathbf{D_{\overline{2}}}$ for {\bf II/10} \newline
iii.e. {\bf II/7} Base-centred lattice, no more restriction to $\Gamma^{0}$, also $\mathbf{D_{\overline{2}}}$ for {\bf II/11} \newline
iii.f. {\bf II/8} Base-primitive, no more restriction to $\Gamma^{0}$, also for {\bf II/12}}.
\label{Fig:8}
\end{figure*}
From (3.8), (3.10) it follows for centred base lattice $\Gamma^0$ that
\begin{equation}
\begin{gathered}
p=0, \ q=1, \ t_1^1=t_1^2 e^{t_3^3}, \ \alpha=0, \ \beta=\frac{1}{t_1^2} (t_3^2 e^{-t_3^3}+t_3^1), \\
\gamma=\frac{2}{t_1^2} \frac{(t_3^1 e^{t_3^3}-t_3^2)}{(e^{2t_3^3}-1)}, \ \ \delta=\frac{2}{t_1^2} \frac{(t_3^2 e^{t_3^3}-t_3^1)}{(e^{2t_3^3}-1)}
\end{gathered} \tag{3.8'}
\end{equation}
From (3.9), (3.10) it follows for primitive base lattice $\Gamma^0$ that
\begin{equation}
\begin{gathered}
N=2p, \ t_2^2=t_1^1\frac{\sqrt{p^2-1}}{q}, \ e^{t_3^3}=p+\sqrt{p^2-1}, \\
\text{i.e.} \ \
\overline{\alpha}=\frac{1+p-\sqrt{p^2-1}}{2t_1^1}\big[t_3^2+t_3^1(p+\sqrt{p^2-1}\big], \\ 
\overline{\beta}=\frac{q(1-p+\sqrt{p^2-1})}{2t_1^1\sqrt{p^2-1}} \big[t_3^2 +t_3^1(p+\sqrt{p^2-1})\big]; \\ 
\overline{\gamma}=\frac{t_3^1+t_3^2}{t_1^1}, \ \  \overline{\delta}=\frac{(t_3^2-t_3^1)q}{t_1^1\sqrt{p^2-1}}. 
\end{gathered} \tag{3.9'}
\end{equation}
\begin{enumerate}
\item[{\bf i.a}] For $t_3^2=t_3^1 e^{2t_3^3}$, we have in (3.8') for centred base lattice $\Gamma^0$ that
\begin{equation}
\begin{gathered}
\alpha=0, \ \beta=\frac{t_3^1}{t_1^2}(e^{t_3^3}+1), \ \gamma=\frac{2t_3^1 e^{t_3^3} (1-e^{t_3^3})}{t_1^2(e^{2t_3^3}-1)}
=\frac{-2t_3^1}{t_1^2} e^{t_3^3} (1+e^{t_3^3}), \\
\delta=\frac{2t_3^1}{t_1^2} \frac{e^{3t_3^3}-1}{e^{2t_3^3}-1}=\frac{2t_3^1 (e^{2t_3^3}+e^{t_3^3}+1)}{t_1^2(e^{t_3^3}+1)},
\text{where} \ \ e^{t_3^3}=\frac{N+\sqrt{N^2-4}}{2}. 
\end{gathered} \notag
\end{equation}
Hence $\beta$ is integer, even  $\beta=1$ can be achived, but then $\gamma=-2e^{t_3^3} \beta$ cannot be integer.
{\it That means}, $t_3^2=t_3^1 e^{2t_3^3}$ {\it realizes} {\bf II/3. i.a.} with $\mathbf{D_{r}}$ in Fig.~8 {\it but not} with $\mathbf{D_{\overline{2}}}$.
\item[{\bf i.b}]
For primitive base lattice $\Gamma^0$ in (3.9'), $t_3^2/t_3^1=e^{2t_3^3}=(p+\sqrt{p^2-1})^2$ yields
\begin{equation}
\begin{gathered}
\overline{\alpha}=\frac{t_3^1}{2t_1^1}(p+\sqrt{p^2-1})[1+p-\sqrt{p^2-1}][1+p+\sqrt{p^2-1}]= \\ =\frac{t_3^1}{t_1^1}(p+\sqrt{p^2-1})(p+1), \\ 
\overline{\beta}=\frac{qt_3^1(p+\sqrt{p^2-1})}{2t_1^1\sqrt{p^2-1}}[1-p+\sqrt{p^2-1}][1+p+\sqrt{p^2-1}]= \\ =\frac{qt_3^1}{t_1^1}(p+\sqrt{p^2-1}), \  
\overline{\gamma}=\frac{t_3^1(e^{2t_3^3}+1)}{t_1^1}=\frac{2pt_3^1}{t_1^1}(p+\sqrt{p^2-1}), \\  
\overline{\delta}=\frac{qt_3^1(e^{2t_3^3}-1)}{t_1^1\sqrt{p^2-1}}=\frac{2qt_3^1}{t_1^1}(p+\sqrt{p^2-1}). 
\end{gathered} \notag
\end{equation}
That means, for $N=2p$, $t_2^2/t_1^1=\frac{\sqrt{p^2-1}}{q}, \ e^{t_3^3}=p+\sqrt{p^2-1}, \ t_3^2/t_3^1=:\overline{x}=(p+\sqrt{p^2-1})^2$,
$t_3^1/t_1^1=:\overline{y}=(p-\sqrt{p^2-1})$, e.g., \ $\overline{\alpha}=p+1, \ \overline{\beta}=q, \ \overline{\gamma}=2p, \ \overline{\delta}=2q$, 
the primitive base lattice $\tau_1(t_1^1,t_1^1,0), \
\tau_2(-t_1^1 \frac{\sqrt{p^2-1}}{q},$ $t_1^1 \frac{\sqrt{p^2-1}}{q},0)$, with $\tau_3(t_1^1 (p-\sqrt{p^2-1}), t_1^1 (p+\sqrt{p^2-1})^2$, \newline
$\log{(1+\sqrt{p^2-1}})$,
realizes the Bravais type {\bf II/4.  i.b} whose point group is automatically $\mathbf{D_{\overline{2}}}$ and has a $z$-sublattice 
$\Gamma_{Oz}$ (Fig.~8.b).
\item[{\bf ii.c}] For $\pm t_3^1 e^{2t_3^3} \ne t_3^2 \ne \mp t_3^1 e^{t_3^3}$ and centred base lattice,  $\beta=1$ can be guaranteed in (3.8'), if 
$t_1^2=t_3^2e^{-t_3^3}+t_3^1$.

Then $t_3^2=t_3^1e^{t_3^3}$ would yield $\gamma=0$, $\delta=\frac{2t_3^1}{t_1^2}=\beta=1$, a convenient realization, just as follows in {\bf iii.e}.

Else in general, the ratio $\overline{x}=t_3^2/t_3^1$ and integer $\beta$ lead in (3.8') to 
\begin{equation}
\begin{gathered}
\alpha=0 \ \ \text{and to} \ \beta=\frac{t_3^1}{t_1^2}(\overline{x}e^{-t_3^3}+1), \\
\gamma=\frac{2(e^{t_3^3}-\overline{x})\beta}{(\overline{x}e^{-t_3^3}+1)(e^{2t_3^3}-1)}, \ \delta=\frac{2(\overline{x}e^{t_3^3}-1)\beta}{(\overline{x}e^{-t_3^3}+1)(e^{2t_3^3}-1)},
\\ \text{with} \ e^{t_3^3}=\frac{N+\sqrt{N^2-4}}{2}, 
\end{gathered} \tag{3.10}
\end{equation}
By eliminating $\overline{x}$ from the above equations
\begin{equation}
\begin{gathered}
\overline{x}=[e^{t_3^3}-\frac{\gamma}{2\beta}(e^{2t_3^3}-1)]\big/(\frac{\gamma}{\beta} \sinh{t_3^3}+1)= \\ =[-1-\frac{\delta}{2\beta}(e^{2t_3^3}-1)]\big/(\frac{\delta}{\beta}
\sinh{t_3^3}-e^{t_3^3})
\end{gathered} \notag
\end{equation}
we get the simple relation for $\beta, \gamma, \delta \ N \in \mathbf{Z}$:
\begin{equation}
\gamma N+2\delta=2\beta \tag{3.11}
\end{equation}
as a criterion for new Bravais type {\bf II/9. ii.c} with point group $\mathbf{D_{\overline{2}}}$ and with $\overline{x}$ above (Fig.~8.c). 
As we look at (3.11) we may have affinely non-equivalent lattices for even $N$ (see Theorem 6). 
At the same time we get realization of type {\bf II/5. ii.c} for other $\overline{x}=t_3^2/t_3^1$ not restricted yet.
\item[{\bf ii.d}] For $\mp t_3^1 e^{t_3^3} \ne t_3^2 \ne \pm t_3^1 e^{2t_3^3}$ and primitive base lattice with $N=2p$,
$t_2^2=\frac{t_1^1\sqrt{p^2-1}}{q}$, $e^{t_3^3}=p+\sqrt{p^2-1}$ in (3.9'), the ratio $\overline{x}=t_3^2/t_3^1$ can guarantee integers $\overline{\alpha}$ and 
$\overline{\beta}$, first. Indeed, we get the equivalence:
\begin{equation}
\begin{gathered}
\frac{t_3^1(\overline{x}+p+\sqrt{p^2-1})}{2t_1^1}=\frac{\overline{\alpha}}{1+p-\sqrt{p^2-1}}= \\ =\frac{\overline{\beta} \sqrt{p^2-1}}{q(1-p+\sqrt{p^2-1})} 
\Leftrightarrow \frac{\overline{\alpha}}{p+1}=\frac{\overline{\beta}}{q}
\end{gathered} \tag{3.12}
\end{equation}
providing realization of type {\bf II/6. ii.d} (Fig.~8.d) with $\overline{x}$ above depending on $t_3^1/t_1^1=:\overline{y}$ yet.

The additional assumptions of (3.9') for integers $\overline{\gamma}$ and $\overline{\delta}$ are necessary for point group $\mathbf{D_{\overline{2}}}$ of our lattice
$\Gamma(\Phi)$. These will fix the ratio $y:=t_3^1/t_1^1$ as follows in the sense of (3.12)
$$\overline{x}=\frac{t_1^1\overline{\alpha}(1+p+\sqrt{p^2-1})}{t_3^1(p+1)}-(p+\sqrt{p^2-1})=\frac{t_1^1 \overline{\gamma}}{t_3^1}-1=\frac{t_1^1}{t_3^1}\frac{\overline{\delta}\sqrt{p^2-1}}{q}+1,$$
\begin{equation}
\begin{gathered}
\text{and so} \ \ \overline{y}=\frac{t_3^1}{t_1^1}= \Big[\overline{\alpha}+\frac{\overline{\alpha}\sqrt{p^2-1}}{p+1}-\overline{\gamma}\Big]\big/(p-1+\sqrt{p^2-1})= \\
=\Big[\overline{\alpha}+\frac{\overline{\alpha}\sqrt{p^2-1}}{p+1}-\frac{\overline{\delta}\sqrt{p^2-1}}{q}\Big]\big/(p-1+\sqrt{p^2-1}), \\ \text{i.e.} \ 
q\Big[\overline{\gamma}-\frac{\overline{\alpha}}{p+1}\Big]=\overline{\delta}(p-1).
\end{gathered} \tag{3.13}
\end{equation}
This determines the Bravais type {\bf II.10. ii.d.} for point group $\mathbf{D_{\overline{2}}}$ (Fig.~8.d) with lattice parameters $\overline{y}$, $\overline{x}$ above.
Again, we may have affinely non-equivalent lattices for given $N$ in this type. 
\item[{\bf iii.e}] For $t_3^2=-t_3^1 e^{t_3^3}$ (similarly $t_3^2= t_3^1 e^{t_3^3}$ for $\delta_{\overline{r}}$) with centred base lattice, formulas (3.8') 
automatically provides
$\alpha=\beta=0$, i.e. point group $\mathbf{D_{r}}$ and Bravais type {\bf II/7. iii.e.}, without restriction by $\tau_3$ to base lattice 
$\Gamma^0$.
\begin{equation}
\begin{gathered}
\gamma=\frac{4t_3^1e^{t_3^3}}{t_1^2(e^{2t_3^3}-1)}, \ \delta=\frac{-2t_3^1(e^{2t_3^3}+1)}{t_1^2(e^{2t_3^3}-1)},
\ \text{with} \ e^{t_3^3}=\frac{N+\sqrt{N^2-4}}{2},\\
\text{lead to} \ \ \frac{\delta}{\gamma}=\frac{N}{2} \ \ \text{and to} \ \ \frac{t_3^1}{t_1^2}=\frac{\gamma}{4}\sqrt{N^2-4}=\frac{\delta}{2N}\sqrt{N^2-4}
\end{gathered} \tag{3.14}
\end{equation}
i.e. to point group $\mathbf{D_{\overline{2}}}$ and Bravais type {\bf II/11. iii.e.}.
\item[{\bf iii.f}] For $t_3^2=-t_3^1 e^{t_3^3}$ (and $t_3^2= t_3^1 e^{t_3^3}$ for $\delta_{\overline{r}}$) and for primitive base lattice 
$\Gamma^0$ in (3.9'), i.e. $N=2p$, $t_2^2=t_1^1 \frac{\sqrt{p^2-1}}{q}$, $e^{t_3^3}=p+\sqrt{p^2-1}$. We get $\overline{\alpha}=0, \ \overline{\beta}=0$ and Bravais type 
{\bf II/8. iii.f.} but also 
\begin{equation}
\begin{gathered}
\overline{\gamma}=-\frac{t_3^1(e^{t_3^3}-1)}{t_1^1}, \ \overline{\delta}=-\frac{t_3^1(e^{t_3^3}+1)q}{t_1^1 \sqrt{p^2-1}},
\ \text{with} \ e^{t_3^3}=p+\sqrt{p^2-1} \\
\text{lead to} \ \ \frac{\overline{\gamma}}{\overline{\delta}}=\frac{p-1}{q} \ \ \text{and to} \ \ \frac{t_3^1}{t_1^1}=\frac{-\overline{\gamma}}{p-1+\sqrt{p^2-1}}
= \\ =\frac{\overline{\gamma}(p-1-\sqrt{p^2-1})}{2(p-1)}.
\end{gathered} \tag{3.15}
\end{equation}
We get the point group $\mathbf{D_{\overline{2}}}$ and Bravais type {\bf II/12. iii.f.}
\end{enumerate}
\medbreak
\item[$\circ$]  {\it The point group is $\mathbf{C_4}$} (Fig.~9. formulas (3.7)

Now $\tau_3^{\gamma_4}\tau_3 \sim (-t_3^2 e^{-t_3^3}+ t_3^1, t_3^1 e^{t_3^3}+t_3^2, 0)\in \Gamma^0$ has to be expressed by basis of 
$\Gamma^0$ $\tau_1(t_1^1,t_1^2,0)$, $\tau_2(t_2^1=-t_1^2,t_2^2=t_1^1,0)$ in integer coordinates $\epsilon, \phi \in \mathbf{Z}$. That means
\begin{equation}
\begin{gathered}
\epsilon t_1^1-\phi t_1^2=-t_3^2 e^{-t_3^3}+t_3^1, \ \ \epsilon t_1^2 + \phi t_1^1= t_3^1 e^{t_3^3}+t_3^2,\ \ \text{i.e.}  \\
\epsilon= \frac{(-t_3^2e^{-t_3^3}+t_3^1)t_1^1+(t_3^1e^{t_3^3}+t_3^2)t_1^2}{t_1^1t_1^1+t_1^2t_1^2}, \ 
\phi=\frac{(t_3^1e^{t_3^3}+t_3^2)t_1^1-(-t_3^2e^{-t_3^3}+t_3^1)t_1^2}{t_1^1t_1^1+t_1^2t_1^2}, \\
\text{and} \ \ \psi t_1^1-\chi t_1^2=2t_3^1, \ \ \psi t_1^2+\chi t_1^1=2t_3^2, \ \ \text{i.e.} \\
\psi=\frac{t_3^1t_1^1+t_3^2t_1^2}{t_1^1t_1^1+t_1^2t_1^2}, \ \ \chi=\frac{t_1^1t_3^2-t_1^2t_3^1}{t_1^1t_1^1+t_1^2t_1^2}, 
\end{gathered} \tag{3.16}
\end{equation}
for $\psi,\chi \in \mathbf{Z}$ by the $\mathbf{C_2}$-action (instead of $\gamma_4^{-1}$ for simplicity).

As before from (2.9) we have for the base lattice
\begin{equation}
\begin{gathered}
-\frac{t_1^2}{t_1^1}= \frac{N-2p-\sqrt{N^2-4}}{2q}, \ \ \frac{t_1^1}{t_1^2}= \frac{N-2p+\sqrt{N^2-4}}{2q}=\frac{-p+e^{t_3^3}}{q}; \\ 
e^{t_3^3}=\frac{N+\sqrt{N^2-4}}{2}, \  \ \text{i.e.} \ \ qq+1=p(N-p).
\end{gathered} \tag{3.17}
\end{equation}
as a criterion, see also (3.7').

By eliminating $\overline{x}:=t_3^2/t_3^1$ and $\overline{y}:=\frac{t_3^1}{t_1^2}$ in (3.16), we have to proceed careful computations by (3.17). We get first
\begin{equation}
\begin{gathered}
\frac{\epsilon}{\overline{y}q}\big[N^2-4+(N-2p)\sqrt{N^2-4}\big]=\overline{\chi}\big[p(N-\sqrt{N^2-4})+ 2(q-1)\big]+ \\ +\big[(N+\sqrt{N^2-4})(1+q)-2p\big], \\
\frac{\phi}{\overline{y}q}\big[N^2-4+(N-2p)\sqrt{N^2-4}\big]=\overline{\chi}\big[2(N-p)+(N-\sqrt{N^2-4})(q-1)\big]+ \\ + \big[(N-p)(N+\sqrt{N^2-4})-2(q+1)\big], \\
\frac{\psi}{\overline{y}q}\big[N^2-4+(N-2p)\sqrt{N^2-4}\big]=\overline{\chi}\big[2q\big]+ \big[N-2p+\sqrt{N^2-4}\big], \\
\frac{\chi}{\overline{y}q}\big[N^2-4+(N-2p)\sqrt{N^2-4}\big]=\overline{\chi}\big[N-2p+\sqrt{N^2-4}\big]+\big[2q\big]. \\
\end{gathered} \tag{3.18}
\end{equation}
Then, say $\phi=0$ can be assumed without lass of generality, since $(t_3^1,t_3^2)$ can be taken $\mod \langle \tau_1,\tau_2 \rangle$. Hence
\begin{equation}
\overline{x}=\frac{2(q+1)-N(N-p)-(N-p)\sqrt{N^2-4}}{N-2p+Nq-(q-1)\sqrt{N^2-4}} \tag{3.19}
\end{equation}
leads to
\begin{equation}
\begin{gathered}
\frac{\epsilon}{2q\overline{y}}=\frac{2}{N-2p+Nq-(q-1)\sqrt{N^2-4}}, \\
\frac{\psi}{2q\overline{y}}=\frac{1-q}{N-2p+Nq-(q-1)\sqrt{N^2-4}}, \\
\frac{\chi}{2q\overline{y}}=\frac{p-N}{N-2p+Nq-(q-1)\sqrt{N^2-4}}, 
\end{gathered} \tag{3.20}
\end{equation}
with integers $\epsilon, \psi, \chi, N, p, q$ as above (3.17), and convenient $\overline{y}$, not uniquely.

This type will be denoted by {\bf II/13} to point group $\mathbf{C_4}$ and called {\it skew tetragonal lattice}. 
\end{enumerate}

We have overviewed the lattices in the main case {\bf II} as well, distinguished 13 Bravais types {\bf II/1-II/13}. All together we have obtained $4+13=17$
Bravais types. $\square$

For summary with realizations e.g. for $N=6$, we formulate
\begin{theorem}
In $\SOL$ space we have 17 Bravais types of lattices. All these are realizable, e.g. for $N=6$ as smallest main parameter with this property. 
Some types of them may have affinely non-equivalent lattice realizations.
\end{theorem}
As a {\it proof} we sketchily give some corresponding lattices by $\tau_1 \sim (t_1^1,t_1^2,0)$, $\tau_2 \sim (t_2^1,t_2^2,0)$, $\tau_3 \sim (t_3^1,t_3^2,t_3^3)$ 
to the proof of Theorem 5, and parameters
$$
N=6, \ p, \ s=N-p, \ q, \ r \in \mathbf{Z} \ \text{with} \ ps-qr=1, \ t_3^3=$$ $$=\log \frac{1}{2} (N+\sqrt{N^2-4})=\log (3+2\sqrt{2}).
$$
{\bf Main case I}; $(t_3^1,t_3^2)=(0,0)$ mod$(\Gamma^0=\langle \tau_1, \tau_2 \rangle)$:
\begin{enumerate}
\item[$\triangleright$] {\it Type {\bf I/1}, point group $\mathbf{C_2}$, \\ Primitive monoclinic lattice} (Fig.~3)
with additional generator the halfturn $\gamma_2$:
\begin{enumerate}
\item[a.] $N=6, \ p=0, \ s=6, \ q=1, \ r=-1,$
\begin{equation}
\begin{gathered}
t_2^1/t_1^1=\frac{1}{2q}(N-2p-\sqrt{N^2-4})=3-2\sqrt{2}, \\ 
t_2^2/t_1^2=\frac{1}{2q}(N-2p+\sqrt{N^2-4})=3+2\sqrt{2}.
\end{gathered} \notag
\end{equation}
\item[b.] $N=6, \ p=1, \ s=N-p=5, \ q=1, \ r=4,$ $t_2^1/t_1^1=2-2\sqrt{2}$, $t_2^2/t_1^2=2+2\sqrt{2}$.
\item[c.] $N=6, \ p=2, \ s=4, \ q=1, \ r=7, \ t_2^1/t_1^1=1-2\sqrt{2}, t_2^2/t_1^2=1+2\sqrt{2}.$

These lattices are affinely equivalent with point group $\mathbf{C_2}$, indeed case b follows from a with
$$
\begin{pmatrix}
u & v \\
w & w 
\end{pmatrix}=\begin{pmatrix}
1 & 0 \\
1 & 1 
\end{pmatrix} \ \ \text{and} \begin{pmatrix}
1 & 0 \\
2 & 1 
\end{pmatrix}
$$
provides the affine equivalence of a and c (fundamental lattice by \cite{AQ}). We shall see that case $N=6, \ p=3, \ s=3, \ q=1, \ r=8$
belongs to Type {\bf I.3}.
\end{enumerate}
\item[$\triangleright$] {\it Type {\bf I/2}, $\mathbf{D_{\overline{2}}}$, \\ Base-centred orthorhombic lattice} (Fig.~4.a):
\begin{enumerate}
\item[a.] $t_2^1=t_1^2, \ t_2^2=t_1^1,\ t_1^1/t_1^2=3+2\sqrt{2}, \ N=6, \ p=0, \ s=N-p=6, \ q=1, \ r=-1$
with additional generators $\delta_r$ and $\gamma_2$, as formerly.
\end{enumerate}
\item[$\triangleright$] {\it Type {\bf I/3}, $\mathbf{D_{\overline{2}}}$, \\ Primitive orthorhombic lattice} (Fig.~4.b):
$0<t_1^2=t_1^1, \ t_2^1=-t_2^2<0, \ t_2^2/t_1^1=2\sqrt{2}/q:$
\begin{enumerate}
\item[a.]  $N=6, \ p=3, \ s=3, \ q=1, \ r=8, \ t_2^2/t_1^1=2\sqrt{2};$
\item[b.]  $N=6, \ p=3, \ s=3, \ q=2, \ r=4, \ t_2^2/t_1^1=\sqrt{2}.$

These lattices are by basis change (2.12), as a simple checking shows.

{\it Thus, they are affinely not equivalent, although they belong to the same Bravais type!}
\end{enumerate}
\item[$\triangleright$] {\it Type {\bf I/4}, $\mathbf{C_4}$, \\ Primitive tetragonal lattice} (Fig.~5):
\begin{enumerate}
\item[a.] $t_2^1=-t_1^2, \ t_2^2=t_1^1,\ t_1^1/t_1^2=\frac{N-2p+\sqrt{N^2-4}}{2q}=1+\sqrt{2}, \ N=6, \ p=1, \ s=5, \ q=r=2.$
\end{enumerate}
\end{enumerate}
\medbreak
{\bf Main case II}; ($t_3^1,t_3^2) \ne (0,0)$ and taken mod$(\Gamma^0=\langle \tau_1, \tau_2 \rangle)$:
\begin{enumerate}
\item[$\triangleright$] {\it Type {\bf II/1}, $\mathbf{Id}$, \\ Primitive triclinic lattice} (Fig.~6):
\begin{enumerate}
\item[a.]  $N=6, \ p=0, \ s=6, \ q=1, \ r=-1, \ t_2^1/t_1^1=3-2\sqrt{2}, \ t_2^2/t_1^2=3+2\sqrt{2};$
\item[b.]  $N=6, \ p=1, \ s=5, \ q=1, \ r=4, \ t_2^1/t_1^1=2-2\sqrt{2}, \ t_2^2/t_1^2=2+2\sqrt{2};$
\item[c.]  $N=6, \ p=2, \ s=4, \ q=1, \ r=7, \ t_2^1/t_1^1=1-2\sqrt{2}, \ t_2^2/t_1^2=1+2\sqrt{2}.$

These lattices are affinely equivalent with point group $\mathbf{Id}$.
\end{enumerate}
\item[$\triangleright$] {\it Type {\bf II/2}, $\mathbf{C_{2}}$, \\ Monoclinic single-side-face-centred lattice} (Fig.~7):
\begin{enumerate}
\item[a.]  $N=6, \ p=0, \ s=6, \ q=1, \ r=-1, \ t_2^1/t_1^1=3-2\sqrt{2}, \ t_2^2/t_1^2=3+2\sqrt{2}, \ t_3^1/t_1^1=t_3^2/t_1^2=\frac{1}{2};$
\item[b.]  $N=6, \ p=1, \ s=5, \ q=1, \ r=4, \ t_2^1/t_1^1=2-2\sqrt{2}, \ t_2^2/t_1^2=2+2\sqrt{2}, \ t_3^1/t_1^1=t_3^2/t_1^2=\frac{1}{2};$
\item[c.]  $N=6, \ p=2, \ s=4, \ q=1, \ r=7, \ t_2^1/t_1^1=1-2\sqrt{2}, \ t_2^2/t_1^2=1+2\sqrt{2}, \ t_3^1/t_1^1=t_3^2/t_1^2=\frac{1}{2},$
can be assumed for additional $\gamma_2$-symmetry to these affinely equivalent (i.e. unimodularly equivalent, see at type {\bf I/1}) lattices, i.e. 
$(\tau_3^{\gamma_2})^{-1}\tau_3=\gamma_2\tau_3=\gamma_2\tau_3^{-1}\gamma_2\tau_3=\tau_1.$
\end{enumerate}
\item[$\triangleright$] {\it Type {\bf II/3. i.a.}, $\mathbf{D_{r}}$, \\ Base-centred with $\Gamma_{Oz}$ lattice} (Fig.~8.a):
\begin{enumerate}
\item[a.]  $t_2^1=t_1^2, \ t_2^2=t_1^1, \ N=6, \ p=0, \ s=6, \ q=1, \ r=-1, \ t_1^1/t_1^2=e^{t_3^3}=\frac{1}{2}(N+\sqrt{N^2-4});
\ \overline{x}:=t_3^2/t_3^1=e^{2t_3^3}$
can be assumed for additional $\delta_r$-symmetry and $\Gamma_{Oz}$ sublattice. That means $(\tau_3^{\delta_r})^{-1}\tau_3 \sim (0,0,2 t_3^3)$,
and for some $\beta \in \mathbf{Z}\setminus \{0\}$ holds $(\tau_3^{\delta_r})\tau_3=\tau_2^{\beta}=
(\tau_1^{\delta_r})^{\beta} \sim (\beta t_1^2, \beta t_1^1,0)$ with 
$\overline{y}:=t_3^1/t_1^2=\beta/(e^{t_3^3}+1)$. Different $\beta$'s yield affinely equivalent lattices.
\end{enumerate}
\item[$\triangleright$] {\it Type {\bf II/4. i.b.}, $\mathbf{D_{\overline{2}}}$, \\ Base-primitive with $\Gamma_{Oz}$ lattice} (Fig.~8.b): 
$t_1^2=t_1^1, \ t_2^1=-t_2^2, \ N=6, \ p=s=3, \ q \ne r$ with $p^2-qr=1, \ t_2^2/t_1^1=\sqrt{p^2-1}/q=2\sqrt{2}/q, \ 
\overline{x}:=t_3^2/t_3^1=e^{2t_3^3}$ moreover $(\tau_3^{\delta_r})\tau_3=\tau_1^{\overline{\alpha}}\tau_2^{\overline{\beta}}, \ 
(\tau_3^{\gamma_2})^{-1}\tau_3= \tau_1^{\overline{\gamma}}\tau_2^{\overline{\delta}}$ and $\overline{y}:=t_3^1/t_1^1=e^{-t_3^3}=
p-\sqrt{p^2-1}=3-2\sqrt{2}$ implies $\overline{\alpha}=p+1=4, \ \overline{\beta}=q, \ \overline{\gamma}=2p=6, \ \overline{\delta}=2q$
for the base lattice:
\begin{enumerate}
\item[a.]  $N=6, \ p=3=s, \ q=1, \ r=8, \ \text{yield} \ \overline{\alpha}=4, \ \overline{\beta}=1, \ \overline{\gamma}=6, \ \overline{\delta}=2;$
\item[b.]  $N=6, \ p=3=s, \ q=2, \ r=4, \ \text{yield} \ \overline{\alpha}=4, \ \overline{\beta}=2, \ \overline{\gamma}=6, \ \overline{\delta}=4,$
\end{enumerate}
i.e. $\overline{\alpha}'=2, \ \overline{\beta}'=1, \ \overline{\gamma}'=3, \ \overline{\delta}'=2$ with $\overline{y}/2=:\overline{y}'$.
These two lattices seemingly are not affine equivalent.
\item[$\triangleright$] {\it Type {\bf II/5. ii.c.}, $\mathbf{D_r}$, \\ Base-centred without $\Gamma_{Oz}$ lattice} (Fig.~8.c): 
$t_2^1=t_1^2, \ t_2^2=t_1^1, \ t_1^1/t_1^2=e^{t_3^3},$ moreover for $\overline{x}:=t_3^2/t_3^1$ hold $\pm e^{2t_3^3}\ne \overline{x} \ne \pm e^{t_3^3}$
and $\alpha=0, \ \beta=\frac{t_3^1}{t_1^2}(\overline{x} e^{-t_3^3}+1)=\overline{y}(\overline{x} e^{-t_3^3}+1)=1$
\begin{enumerate}
\item[a.]  $N=6, \ p=0, \ s=6, \ q=1, \ r=-1, \ \text{and} \ \gamma$ or $\delta$ is not integer in (3.10) for given $\overline{x}$ above.
Then $\mathbf{D_r}$ is the point group (not $\mathbf{D_{\overline{2}}}$).
\end{enumerate}
\item[$\triangleright$] {\it Type {\bf II/9. ii.c.}, $\mathbf{D_{\overline{2}}}$, \\ Base-centred without $\Gamma_{Oz}$ lattice} (Fig.~8.c):
As above, but $\overline{x}:=t_3^2/t_3^1$ and $\overline{y}:=t_3^1/t_1^2$ with $\beta=1$ satisfy (3.10) with integer $\gamma$ and $\delta$ with
$\gamma N+2 \delta =2$ \ (3.11).
\begin{enumerate}
\item[a.]  As above at {\bf II/5}, but for some $\overline{x}<0, \ \gamma=1, \ \delta=-2; \ \beta=1$ with a convenient $\overline{y}$ provide a lattice with
$\mathbf{D_{\overline{2}}}$-symmetry. 
\item[b.] As above, but for some $\overline{x}>0, \ \gamma=-1, \ \delta=4; \ \beta=1$ with a convenient $\overline{y}$ provide another lattice, again
with $\mathbf{D_{\overline{2}}}$-symmetry. 
\end{enumerate}
\item[$\triangleright$] {\it Type {\bf II/6. ii.d.}, $\mathbf{D_{r}}$, \\ Base-primitive without $\Gamma_{Oz}$ lattice} (Fig.~8.d): 
$0<t_1^2=t_1^1, \ t_2^1=-t_2^2<0, \ t_2^2/t_1^1=\sqrt{p^2-1}/q,$ moreover for  \ 
$\overline{x}:=t_3^2/t_3^1$ hold $\pm e^{2t_3^3} \ne \overline{x} \ne \pm e^{t_3^3}$ in equations (3.12) with $\overline{y}:=t_3^1/t_1^1$ 
and intger $\overline{\alpha}$ and $\overline{\beta}$. But $\overline{\gamma}$ or $\overline{\delta}$ is not integer in (3.9'), (3.13):
\begin{enumerate}
\item[a.]  $N=6, \ p=3=s, \ q=1, \ r=8, \ \overline{\alpha}=p+1=4, \ \overline{\beta}=q=1;$
\item[b.]  $N=6, \ p=3=s, \ q=2, \ r=4, \ \overline{\alpha}=\frac{p+1}{2}=2, \ \overline{\beta}=\frac{q}{2}=1.$
\end{enumerate}
\item[$\triangleright$] {\it Type {\bf II/10. ii.d.}, $\mathbf{D_{\overline{2}}}$, \\ Base-primitive without $\Gamma_{Oz}$ lattice} (Fig.~8.d): 
As previously at {\bf II/6}, but 
$\overline{x}:=t_3^2/t_3^1$ and $\overline{y}:=t_3^1/t_1^1$ with intger $\overline{\alpha}$ and $\overline{\beta}$, $\frac{\overline{\alpha}}{p+1}=
\frac{\overline{\beta}}{q}$, satisfy also (3.13) with integer $\overline{\gamma}$ and $\overline{\delta}, 
q\big(\overline{\gamma}-\frac{\overline{\alpha}}{p+1}\big)=\overline{\delta}(p-1):$ 
\begin{enumerate}
\item[a.]  $N=6, \ p=3=s, \ q=1, \ r=8, \ \overline{\alpha}=p+1=4, \ \overline{\beta}=q=1, \ \overline{\gamma}=3, \ \overline{\delta}=1$;
\item[b.]  $N=6, \ p=3=s, \ q=2, \ r=4, \ \overline{\alpha}=4, \ \overline{\beta}=2, \ \overline{\gamma}=2, \ \overline{\delta}=1$.
\end{enumerate}
Of course, these are in both cases with convenient $\overline{x}, \overline{y}$, respectively, by (3.13).
\item[$\triangleright$] {\it Type {\bf II/7. iii.e.}, $\mathbf{D_{r}}$, \\ Base-centred, no more restriction to $\Gamma^0$ lattice} (Fig.~8.e): 
$t_2^1=t_1^2, \ t_2^2=t_1^1, \ t_1^1/t_1^2=e^{t_3^3},$ moreover 
$\overline{x}:=t_3^2/t_3^1=-e^{t_3^3}, \ \ \alpha=\beta=0$ automatically holds:
\begin{enumerate}
\item[a.]  $N=6, \ p=0, \ N=s=6, \ q=1, \ r=-1,$ with  $\overline{x}=-e^{t_3^3}, \ \overline{y}=t_3^1/t_1^2$ is not restricted yet.
\end{enumerate}
\item[$\triangleright$] {\it Type {\bf II/11. iii.e.}, $\mathbf{D_{\overline{2}}}$, \\ Base-centred, no more restriction to $\Gamma^0$ lattice} (Fig.~8.e): 
As above at {\bf II/7}, but $\overline{y}:=\frac{t_3^1}{t_1^1}=\frac{\gamma}{4}\sqrt{N^2-4}=\frac{\delta}{2N}\sqrt{N^2-4}$ 
with intger $\gamma$ and $\delta$, guarantees the point group $\mathbf{D_{\overline{2}}}$:
\begin{enumerate}
\item[a.]  $N=6, \ p=0, \ s=6, \ q=1, \ r=-1, \ \gamma=1, \ \delta=3$ with $\overline{y}=\frac{t_3^1}{t_1^2}=\sqrt{2}$.
\end{enumerate}
\begin{figure*}[ht]
\includegraphics[width=6cm]{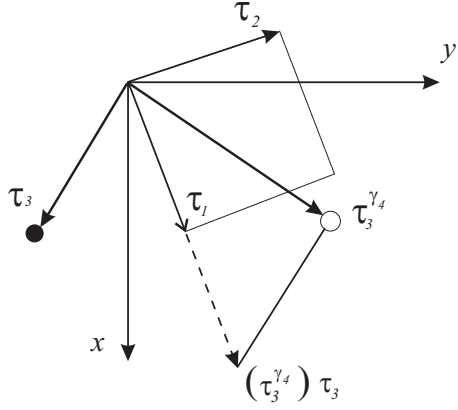} 
\caption{Symbolic picture of the skew tetragonal lattice to point group $\mathbf{C_{4}}$}.
\label{Fig:8}
\end{figure*}
\item[$\triangleright$] {\it Type {\bf II/8. iii.f.}, $\mathbf{D_{r}}$, \\ Base-primitive, no more restriction to $\Gamma^0$ lattice} (Fig.~8.f):
$0<t_1^2=t_1^1, \ t_2^1=-t_2^2<0, \ t_2^2/t_1^1=\sqrt{p^2-1}/q,$ moreover for  \ 
$\overline{x}:=t_3^2/t_3^1= -e^{2t_3^3}$ alredy yields $\mathbf{D_{r}}$, $\overline{\alpha}=0=\overline{\beta}$ without restriction to
$\overline{y}:= t_3^1/t_1^1$:
\begin{enumerate}
\item[a.]  $N=6, \ p=3=s, \ q=1, \ r=8, \text{with} \ \overline{x}=t_3^2/t_3^1=-e^{t_3^3};$
\item[b.]  $N=6, \ p=3=s, \ q=2, \ r=4, \ \overline{x}=-e^{t_3^3}.$
\end{enumerate}
\item[$\triangleright$] {\it Type {\bf II/12. iii.f.}, $\mathbf{D_{\overline{2}}}$, \\ Base-primitive, 
no more restriction to $\Gamma^0$ lattice} (Fig.~8.f):
As above at {\bf II/8}, but 
$\overline{y}:=t_3^1/t_1^1=\overline{\gamma}(p-1-\sqrt{p^2-1})\big/(2(p-1))$ with  $\frac{\overline{\gamma}}{\overline{\delta}}=
\frac{p-1}{q}$, for integer $\overline{\gamma}, \overline{\delta}:$ 
\begin{enumerate}
\item[a.]  $N=6, \ p=3=s, \ q=1, \ r=8, \ \overline{\gamma}=2, \ \overline{\delta}=1, \ \overline{y}=2(2-2\sqrt{2})$;
\item[b.]  $N=6, \ p=3=s, \ q=2, \ r=4, \ \overline{\gamma}=1, \ \overline{\delta}=1, \ \overline{y}=2-2\sqrt{2}$.
\end{enumerate}
\item[$\triangleright$] {\it Type {\bf II/13.}, $\mathbf{C_{4}}$, \\ Skew tetragonal lattice} (Fig.~9):
$t_2^1=-t_1^2, \ t_2^2=t_1^1, \ t_1^1/t_1^2=\frac{-p+e^{t_3^3}}{q}, \ e^{t_3^3}=\frac{1}{2}(N+\sqrt{N^2-4}), \ qq+1=p(N-p):$
\begin{enumerate}
\item[a.]  $N=6, \ p=1, \ s=5, \ q=r=2, \ \overline{x}=t_3^2/t_3^1=\frac{6-30-5\sqrt{32}}{4+12-\sqrt{32}}$ by (3.19),
$\phi=0, \ \epsilon=2, \ \psi=-1, \ \chi=-5$ by (3.20), with \ $\overline{y}:=\frac{t_3^1}{t_1^2}=\frac{1}{4}(4+12-\sqrt{32})=4-\sqrt{2}. \ \square$
\end{enumerate}
\end{enumerate}

{\bf Closing remarks:}

The above classification of $\SOL$-lattices into Bravais types may open certain investigations in discrete geometry 
(see e.g. \cite{Sz10}) and topology (see \cite{CMSS09}). The number theoretical problems for finer affine classification,
sketched in this paper, will be left for future works.

Our next intention is to find a compact $\SOL$ manifold with 2-generator fundamental group, different from one generated by 2 translations $\tau_1$ and 
$\tau_3$ above.

\end{document}